\input amstex
\documentstyle{amsppt}
\magnification=1200
\hoffset=-0.5pc
\nologo
\vsize=57.2truepc
\hsize=38.5truepc

\spaceskip=.5em plus.25em minus.20em

\define\hpsi{\psi}
\define\Bobb{\Bbb}
\define\fra{\frak}
\define\KK{H}
\define\kk{h}
 \define\armcusgo{1}
 \define\brionone{2}
 \define\drivhall{3}
 \define\emmroeme{4}
 \define\gotayone{5}
 \define\gotaythr{6}
 \define\gotaysix{7}
 \define\guistetw{8}
 \define\guisteth{9}
\define\bhallone{10}
\define\hartsboo{11}
 \define\howeone{12}
\define\poiscoho{13}
 \define\souriau{14}
    \define\srni{15}
\define\claustha{16}
\define\oberwork{17}
 \define\kaehler{18}
 \define\severi{19}
 \define\descent{20}
\define\kirwaboo{21}
\define\kostaone{22}
\define\kralyvin{23}
\define\meinrtwo{24}
\define\naramtwo{25}
 \define\nessone{26}
\define\ramadthr{27}
\define\roberone{28}
\define\schlione{29}
\define\sjamatwo{30}
\define\sjamafou{31}
\define\sjamlerm{32}
\define\sniabook{33}
\define\sniatone{34}
\define\sniawein{35}
\define\ctelethr{36}
\define\woodhous{37}
\define\zheloboo{38}
\topmatter
\title K\"ahler quantization and reduction
\endtitle
\author Johannes Huebschmann
\endauthor
\affil 
Universit\'e des Sciences et Technologies
de Lille
\\
U. F. R. de Math\'ematiques
\\
CNRS-UMR 8524
\\
F-59 655 VILLENEUVE D'ASCQ, France
\\
Johannes.Huebschmann\@math.univ-lille1.fr
\endaffil
\date{April 5, 2004}
\enddate
\abstract{
Exploiting a notion of K\"ahler structure on a stratified space introduced 
elsewhere we show that, in the K\"ahler case, reduction after quantization 
coincides with quantization after reduction: Key tools developed for that 
purpose are stratified polarizations and stratified prequantum modules, the 
latter generalizing prequantum bundles. These notions encapsulate, in 
particular, the behaviour of a polarization and that of a prequantum bundle 
across the strata. Our main result says that, for a positive K\"ahler manifold 
with a hamiltonian action of a compact Lie group, when suitable additional 
conditions are imposed, reduction after quantization coincides with 
quantization after reduction in the sense that not only the reduced and 
unreduced quantum phase spaces correspond but the (invariant) unreduced and 
reduced quantum observables as well. Over a stratified space, the appropriate 
quantum phase space is a costratified Hilbert space in such a way that the 
costratified structure reflects the stratification. Examples of stratified 
K\"ahler spaces arise from the closures of holomorphic nilpotent orbits 
including angular momentum zero reduced spaces and from representations of 
compact Lie groups. For illustration, we carry out K\"ahler quantization on 
various spaces of that kind including singular Fock spaces.}
\endabstract

\address{\smallskip
\noindent
USTL, UFR de Math\'ematiques, CNRS-UMR 8524
\newline\noindent
F-59 655 Villeneuve d'Ascq C\'edex,
France
\newline\noindent
Johannes.Huebschmann\@math.univ-lille1.fr}
\endaddress
\subjclass
\nofrills{{\rm 2000}
{\it Mathematics Subject Classification}.\usualspace}
{
14L24 
14L30 
17B63 
17B65 
17B66 
17B81 
32C20 
32Q15 
32S05 
32S60 
53D17
53D20
53D50
81S10 
}
\endsubjclass
\keywords{Poisson manifold, Poisson algebra,
Poisson cohomology, holomorphic quantization,
reduction and quantization,
stratified K\"ahler space,
geometric quantization,
quantization on a space with singularities,
normal complex analytic space,
constrained system,
invariant theory,
real Lie algebra of hermitian type}
\endkeywords

\endtopmatter
\document
\leftheadtext{Johannes Huebschmann}
\rightheadtext{K\"ahler quantization and reduction}

\medskip
\noindent
{\bf Introduction}
\smallskip\noindent
The relationship between unitary representations of a compact Lie group
$G$ and K\"ahler quantization on
smooth compact hamiltonian $G$-spaces
has received much attention.
In this paper, we will develop a similar theory
for hamiltonian $G$-spaces which are not necessarily
smooth manifolds. Our motivation comes from physics:
Given a quantizable system with constraints,
the question arises whether reduction after 
quantization coincides with quantization after reduction, so that it 
would then make no difference whether reduction is imposed before or 
after quantization.
This question goes back to the early days of quantum mechanics
and appears already in  {\smc Dirac's}
work on the electron and positron.
In the present paper we will 
show that, indeed, in a suitable sense,
in the framework of K\"ahler quantization,
reduction after quantization is equivalent to
quantization after reduction.
\smallskip
When  the unreduced phase space is a quantizable
smooth symplectic manifold and when the symmetries
can be quantized as well
reduction after quantization 
is well defined
within the usual framework of
geometric quantization.
Up to now,
the available methods 
have been insufficient to attack the problem of
quantization of reduced observables, though,
once the reduced phase space is no longer a smooth manifold;
we will refer to this situation as the {\it singular case\/}.
The singular case is the rule rather than the exception.
For example, simple classical mechanical systems and
the solution spaces of classical field theories
involve singularities;
see e.~g. \cite\armcusgo\ and the references there.
In the presence of singularities, restricting quantization
to a smooth open dense part, the \lq\lq top stratum\rq\rq, can result in a 
loss of information and may in fact lead to inconsistent results, cf. (4.12) 
below. This kind of phenomenon has been long known, perhaps
as a folk-lore observation,
but it does not seem easy to trace it down explicitly in the literature;
cf. e.~g. the discussion in \cite\emmroeme\ and \cite\gotayone.
To overcome these difficulties on the classical level,
in a predecessor to this paper \cite\kaehler,
we isolated a certain class
of \lq\lq K\"ahler spaces with singularities\rq\rq,
which we call stratified K\"ahler spaces.
On such a space,
the complex analytic structure alone is unsatisfactory
for issues related with
quantization because it overlooks the requisite Poisson structures.
In this paper, we generalize ordinary K\"ahler quantization 
to a quantization scheme over
(complex analytic) stratified K\"ahler spaces.
\smallskip
We now explain briefly and informally our approach:
Consider 
a stratified symplectic space $(N,C^{\infty}(N),\{\cdot,\cdot\})$.
The Poisson structure 
encapsulates the mutual positions of the symplectic structures
on the strata of $N$.
Likewise,  a {\it stratified K\"ahler polarization\/}
(cf.  \cite\kaehler)
induces ordinary K\"ahler polarizations on the strata and
{\it encapsulates the mutual positions of these polarizations
on the strata\/}.
A complex polarization can {\it no longer
be thought of as being given by the $(0,1)$-vectors
of a complex structure\/}, though;
see Remark 3.7.3 in \cite\kaehler.
A suitable notion of
{\it prequantization\/}, phrased in terms {\it prequantum modules\/} 
introduced in \cite\souriau,
yields the requisite representation of the 
Poisson algebra; in particular, this representation satisfies the
Dirac condition. 
The concept of
stratified K\"ahler polarization then takes care of the 
irreducibility problem, as  does an ordinary polarization
in the smooth case.
Over a stratified space,
the appropriate quantum phase space is what we call a
{\it costratified\/} Hilbert space; this is a system of Hilbert spaces,
one for each stratum, which arises from quantization
on the closure of that stratum, the stratification
provides linear maps between these Hilbert spaces
reversing the partial ordering among the strata,
and these linear maps are compatible with the quantizations. 
After these preparations we show
in Section 3 below
that, for a positive K\"ahler manifold
with a hamiltonian action of a compact Lie group,
when suitable additional conditions are imposed, reduction after 
quantization coincides with quantization after reduction in the sense that not 
only the invariant unreduced 
and reduced 
quantum phase spaces correspond but the 
invariant unreduced and reduced quantum observables as well.
The correspondence between the
invariant unreduced and reduced 
quantum phase spaces is given by a suitable map
which we will refer to as a {\it state space comparison map\/};
in the literature, the phrase \lq quantization commutes with reduction\rq\ 
is usually interpreted as the requirement that 
a state space comparison map
be an isomorphism (of complex vector spaces).
The comparison of invariant unreduced
and reduced quantizable classical observables
we will give below (Theorem 3.6) does not rely
on the state space
comparison map 
(written as $\rho$ in Section 3 below)
being an isomorphism, though,
that is, 
the phrase \lq quantization commutes with reduction\rq\ 
may be given a consistent meaning
whether or not the state space
comparison map is an isomorphism,
and the question arises
whether there are
reduced states which do {\it not\/} arise
from (invariant) unreduced states.
See Remark 3.9 below.
\smallskip
We illustrate our approach with a number of examples
involving what we call {\it singular Fock spaces\/}.
See Section 4 below for details.
In particular, exploiting
a suitable notion of momentum mapping
which is Poisson even when defined on a space with singularities,
we show how
the relationship between unitary representations of a compact Lie group
and K\"ahler quantization 
extends to  certain singular cases.
Momentum mappings defined on not necessarily
smooth spaces occur already in the literature,
for example in \cite\brionone\ 
and \cite\nessone\ (at least implicitly),
but these momentum mappings
do not involve global Poisson structures
encapsulating the singular behaviour.
See Remark 4.14 below.
Issues related with the metaplectic correction in singular situations
will be addressed elsewhere.
\smallskip
In recent years considerable work has been devoted to the 
comparison of the reduced and unreduced
quantum phase 
spaces and to variants thereof,
see \cite\sjamafou\ 
for an overview and the literature there, 
as well as e.~g. \cite\naramtwo\  and \cite\ramadthr, \cite\ctelethr\ 
for generalizations to higher dimensional sheaf cohomology.
The actual quantization of corresponding
classical observables has received much less attention, though;
cf. e.~g. \cite\gotaythr\ and \cite\sniatone\ 
where the problem has been studied for hamiltonain $G$-spaces
endowed with a $G$-invariant real polarization.
\smallskip
This paper is part of a program aimed at
developing a satisfactory quantization procedure on certain moduli
spaces including spaces of possibly twisted representations
of the fundamental group of a surface in a compact Lie group.
In \cite\kaehler\  we have shown  that these spaces are indeed
normal K\"ahler spaces;
see our expository papers \cite{\srni--\oberwork} and the
literature there for the stratified symplectic structure.
In physics language, the relevant 
quantum phase spaces are spaces of conformal blocks.
\smallskip
I am indebted to A. Weinstein for discussions,
and for his encouragement to carry out the research program
the present paper is part of. 
I am grateful to T. R. Ramadas for having pointed out to me
some relevant literature and to M. Schmidt 
and J. \'Sniatycki for a number of questions
which helped improve the exposition.

\medskip\noindent {\bf 1. Prequantization on spaces with singularities}
\smallskip\noindent
To develop
prequantization over stratified symplectic spaces and
to describe the behaviour of prequantization under reduction, 
we will introduce {\it stratified prequantum 
modules\/} over stratified symplectic spaces.
A stratified
prequantum 
module
determines what we call a {\it costratified prequantum space\/}
but the two notions, though closely related, should not be confused.
For intelligibility, we reproduce first the concept of
prequantum bundle in a language tailored to our purposes.
\smallskip\noindent
(1.1) {\sl Prequantum bundles\/.}
Let $(N,\sigma)$ be a quantizable symplectic manifold,
let
$(C^{\infty}(N),\{\cdot,\cdot\})$
be its symplectic Poisson algebra,
and let  $\zeta\colon \Lambda \to N$ be a prequantum bundle
for $(N,\sigma)$.
Thus, when
the operator of covariant derivative for the requisite connection
is written as
$\nabla
\colon \Omega^0(N,\zeta)
@>>>
\Omega^1(N,\zeta)$,
the curvature $K_\nabla$ coincides with
$-i \sigma$.
Here the convention is that
the curvature $K_\nabla$ of the connection
$\nabla$
is characterized by the formula
$[\nabla_X,\nabla_Y]= \nabla_{[X,Y]} + K_\nabla(X,Y)$
where $X$ and $Y$ are arbitrary smooth vector fields on $N$
and,
for  
a smooth vector field
$X$ on $N$, $\nabla_X$ is the operator which
assigns the smooth complex valued section
$\nabla_{X} s = (\nabla (s))(X)$
to 
a smooth complex valued section $s$ of $\zeta$.
Henceforth we will often write
the space $\Omega^0(N,\zeta)$
of {\it smooth\/} complex sections of $\zeta$
as $\Gamma^{\infty}(\zeta)$. 
Ordinary prequantization proceeds by means of Kostant's formula
$$
\widehat f (s) = - i \nabla_{\{f,\cdot\}} (s) + fs,
\quad
f \in C^{\infty}(N),\ s \in \Gamma^{\infty}(\zeta) .
\tag1.1.1
$$
(We write $\nabla_{\{f,\cdot\}}$ rather than a corresponding expression
involving the hamiltonian vector field $X_f$ of 
the function $f$ since, in accordance 
with Hamilton's equations, the Hamiltonian vector field of $f$ is given by 
the operator ${\{\cdot\,,f\}}$.)
Associating
$\widehat f$ to $f$
yields a representation
on 
$\Gamma^{\infty}(\zeta)$
of the  real Lie algebra
which underlies $C^{\infty}(N)$;
here $\Gamma^{\infty}(\zeta)$
is viewed as a complex vector space,
and the elements of 
$C^{\infty}(N)$ are represented
by $\Bobb C$-linear transformations
so that the constants in $C^{\infty}(N)$ 
act by {\it multiplication\/}
and the {\it Dirac condition holds\/} (see (1.2.6) and (1.2.7) below).
The physical constant
$\hbar$ 
is absorbed
in the symplectic or, what amounts to the same, Poisson structure;
see Remark 1.2.10 below.
\smallskip\noindent
(1.2) {\sl Prequantum modules\/.}
Let $(A,\{\cdot,\cdot\})$ be an arbitrary real Poisson algebra.
Recall from \cite\poiscoho\ and \cite\souriau\ that
the Poisson structure
$\{\cdot,\cdot\}$ determines
an $(\Bobb R,A)$-Lie algebra structure
$([\cdot,\cdot],\pi^{\sharp})$
on 
the $A$-module $D_A$ of formal differentials for $A$.
Here 
$\pi =\pi_{\{\cdot,\cdot\}}\colon D_A \otimes D_A \to A$ is the 2-form given by
$\pi(da,db) = \{a,b\}$ ($a,b \in A$), the morphism
$\pi^{\sharp}$ from $D_A$ to $\roman{Der}(A) = \roman{Hom}_A(D_A,A)$
is the adjoint of $\pi$, and
the bracket $[\cdot,\cdot]$ on $D_A$ is given by the formula
$$
[adu,bdv] = a \{u,b\}dv + b \{a,v\}du + ab d\{u,v\},
\quad a,b,u,v \in A.
$$
See \cite\poiscoho\ and \cite\souriau\
for details.
We write the resulting 
$(\Bobb R,A)$-Lie algebra as
$D_{\{\cdot,\cdot\}}$;
thus the pair $(A,D_{\{\cdot,\cdot\}})$ is a Lie-Rinehart algebra.
For intelligibility we recall that,
given a Lie-Rinehart algebra $(A,L)$, 
the Lie algebra $L$ together with the additional structure
is referred to as an $(R,A)$-{\it Lie algebra\/}.
\smallskip
The Poisson 2-form $\pi_{\{\cdot,\cdot\}}$
determines an {\it extension\/}
$$
0 
@>>> 
A
@>>> 
\overline L^{\roman{a}}_{\{\cdot,\cdot\}}
@>>> 
D_{\{\cdot,\cdot\}}
@>>> 
0
\tag1.2.1
$$
of $(\Bobb R,A)$-Lie algebras which is central 
as an extension of ordinary Lie algebras;
in particular, on the kernel $A$, the Lie bracket is trivial.
At the risk of making a mountain out of a molehill we note that
here and below,
according to standard
conventions,
\lq\lq
\,$
0 
@>>> 
U
@>>> 
V
$\,\rq\rq\
and
\lq\lq
\,$
V 
@>>>
W
@>>>
0
$\, \rq\rq\ 
signify that 
$U@>>> V$
and  
$V@>>> W$
are injective and surjective
morphisms of $A$-modules, respectively;
in particular, $0$ is {\it not\/} considered as an $(\Bobb R,A)$-Lie algebra.
Moreover, as $A$-modules,
$$
\overline L^{\roman{a}}_{\{\cdot,\cdot\}} = A \oplus D_{\{\cdot,\cdot\}},
\tag1.2.2
$$
and the Lie bracket on $\overline L^{\roman{a}}_{\{\cdot,\cdot\}}$ 
is given by
$$
[(a,du),(b,dv)] =
\left(
\{u,b\}+ \{a,v\} - \{u,v\}, d\{u,v\}
\right) ,\quad
a,b,u,v \in A.
\tag1.2.3
$$
The superscript \lq\lq \,${\ldots}^{\roman{a}}$\,\rq\rq\ 
is intended to refer to \lq\lq algebraic\rq\rq\ (this superscript does not
occur in \cite\poiscoho\ and \cite\souriau), 
and we have written \lq\lq $\overline L^{\roman{a}}$\rq\rq\ 
rather than simply $L^{\roman{a}}$ to indicate that
the extension (1.2.1) represents the {\it negative\/} of the class of
$\pi_{\{\cdot,\cdot\}}$
in Poisson cohomology
$\roman H_{\roman{Poisson}}^2(A,A)$, cf. \cite\poiscoho.
When $(A,\{\cdot,\cdot\})$ is the smooth symplectic Poisson algebra
of an ordinary smooth symplectic manifold, cf. (1.1),
(perhaps) up to sign, the class of $\pi_{\{\cdot,\cdot\}}$
comes down to the cohomology class represented by the symplectic structure.
Plainly the extension (1.2.1) determines an extension of Lie-Rinehart algebras
as well, the algebra variable being fixed.
\smallskip
Extending terminology introduced in \cite\souriau,
given an
$(A\otimes \Bobb C)$-module
$M$,
we refer to
an 
$(A,\overline L^{\roman{a}}_{\{\cdot,\cdot\}})$-module
structure 
$$
\chi
\colon 
\overline L^{\roman{a}}_{\{\cdot,\cdot\}}
\longrightarrow
\roman{End}_{\Bobb R}(M)
\tag1.2.4
$$
on $M$  
as an 
{\it algebraic prequantum module structure for\/}
$(A,{\{\cdot,\cdot\}})$
provided 
(i) the values of $\chi$ lie in
$\roman{End}_{\Bobb C}(M)$,
that is to say, 
for $a \in A$ and $\alpha \in 
D_{\{\cdot,\cdot\}}$,
the operators $\chi(a,\alpha)$ are complex linear
transformations,
and (ii)
for every $a\in A$, with reference to the decomposition (1.2.2), we have
$$
\chi(a,0) = i\,a\,\roman{Id}_M.
\tag1.2.5
$$
In \cite\poiscoho, the terminology \lq prequantum module structure\rq\ 
is used for what we here call algebraic prequantum module structure.
A pair $(M,\chi)$
consisting of 
an $(A\otimes \Bobb C)$-module $M$ and an algebraic prequantum module structure
will henceforth be referred to as an {\it algebraic prequantum module\/}
(for $(A,\{\cdot,\cdot\})$.
\smallskip
Prequantization now proceeds 
in the following fashion,
cf. \cite\poiscoho:
The assignment to $a \in A$ of
$(a,da) \in
\overline L^{\roman{a}}_{\{\cdot,\cdot\}}$
yields a morphism $\iota$ of real Lie algebras
from
$A$ to
$\overline L^{\roman{a}}_{\{\cdot,\cdot\}}$;
thus, for any algebraic prequantum module $(M,\chi)$,
the composite of $\iota$ with $-i \chi$
is a representation
$a \mapsto \widehat a$
of the $A$ underlying real Lie algebra
having $M$, viewed as a complex vector space,
as its representation space;
this is a representation by $\Bobb C$-linear operators 
so that any constant acts by multiplication,
that is, 
for any real number $r$,
viewed as a member of $A$,
$$
\widehat r = r \,\roman{Id}
\tag1.2.6
$$
and so that, for $a,b \in A$,
$$
\widehat {\{a,b\}} = i\,[\widehat a,\widehat b]
\qquad
\text{(the Dirac condition).}
\tag1.2.7
$$
More explicitly, these operators are given by the formula
$$
\widehat a (x) = \frac 1 i \chi(0,da) (x) + ax,
\quad
a \in A,\ x \in M,
\tag1.2.8
$$
which we shall henceforth refer to as the 
{\it prequantization formula\/}.
\smallskip
The 
interpretation of 
quantum mechanics requires
observables
to be represented by symmetric operators
(after introduction of a suitable Hilbert space structure),
and this forces the factor $i$ in the Dirac condition 
(1.2.7) (since the ordinary commutator of two symmetric operators
is skew); this factor $i$, in turn, forces 
multiplication of the structure map $\chi$ of a prequantum module
by $-i$.
\smallskip
Symmetries are to be quantized
by skew symmetric operators, though, that is,
when a classical observable 
$a \in A$
is viewed as an {\it infinitesimal symmetry\/},
the corresponding infinitesimal quantum symmetry
is given by the operator
$\widetilde a = i \widehat a$, that is,
$$
\widetilde a (x) =  (\chi(a,da)) (x) = (\chi(0,da)) (x) + iax,
\quad
a \in A,\ x \in M.
\tag1.2.9
$$
Thus, when $\widehat a$ is self-adjoint (with reference to an appropriate
Hilbert space structure, perhaps on
a suitable subspace of $M$), it generates a 1-parameter group of unitary 
transformations.

\smallskip\noindent
{\smc Remark 1.2.10.}
In the situation of (1.1), 
let $\omega = \hbar \sigma$, so that
$\sigma = \frac{\omega}{\hbar}$ ($=\frac{2 \pi \omega}h$).
Then,
using superscripts to indicate with reference to which 
symplectic structure hamiltonian vector fields and Poisson brackets
are to be taken, given functions $f$ and $g$, we have
$\{f,g\}^{\sigma} = \hbar \{f,g\}^{\omega}$,
the prequantization formula (1.1.1)
may be written as
$$
\widehat f (s) = - i \hbar\nabla_{\{f,\cdot\}^{\omega}} (s) + fs,
\tag1.2.11
$$
and the formula (1.1.5) is equivalent to 
$$
i [\widehat a,\widehat b] = \hbar\widehat {\{a,b\}^{\omega}} 
\tag1.2.12
$$
which is more common in the physics literature.
Now
the quantizability of $\sigma$, i.~e.
the existence of a prequantum bundle $\zeta$,
is equivalent to the requirement that, under the integration map
$\roman H^2_{\roman{de Rham}}(N,\Bobb R) \to
\roman H^2_{\roman{singular}}(N,\Bobb R)$,
the class of $\sigma$ go to $2\pi$ times the image of an integral class,
and we may then think of the class of
$\frac{\omega}h$ as this integral class.
However, the original phase space Poisson bracket is that coming from 
$\sigma$ and {\it not\/} from $\omega$ nor $\frac \omega h$.
We therefore believe that our formula (1.2.9) is more appropriate than
corresponding ones in the literature
involving a factor $2 \pi$ (and hence the Poisson bracket
coming from $\frac \omega h$).
\smallskip\noindent
(1.3) {\sl Stratified symplectic spaces\/.}
Let $N$ be a stratified symplectic space, and let
$(A,\{\cdot,\cdot\})$
be its stratified symplectic Poisson algebra; 
a special case would
be
the ordinary symplectic Poisson algebra of 
a smooth symplectic manifold.
Dividing out the formal differentials
in $D_A$
that {\it vanish at every point of $N$\/}
yields 
the $A$-module
$\Omega^1(N)$
which serves as a module of differentials for $A$ as well \cite\kralyvin.
Here a formal differential $\alpha$ in $D_A$
{\it vanishes at the point\/} $w$ of $N$ provided
$\alpha$ goes to zero under the epimorphism
$D_A \to \Bobb R \otimes_AD_A$ induced by the point $w$.
Equivalently, $\Omega^1(N)$ is the quotient 
of $D_A$ by the formal differentials $\alpha$ having the property that
$(X,\alpha)$ is zero for {\it every\/} derivation $X$ of $A$.
The
$(\Bobb R,A)$-Lie algebra
structure 
$([\cdot,\cdot],\pi^{\sharp})$
on $D_A$
descends to
an $(\Bobb R,A)$-Lie algebra structure on
$\Omega^1(N)$;
we write the resulting
$(\Bobb R,A)$-Lie algebra
as
$\Omega^1(N)_{\{\cdot,\cdot\}}$.
Thus the 
$(\Bobb R,A)$-Lie algebra
$\Omega^1(N)_{\{\cdot,\cdot\}}$
consists of
the $A$-module $\Omega^1(N)$
endowed with the induced bracket
$[\cdot,\cdot]$ and structure map $\pi^{\sharp}$
from $\Omega^1(N)$ to $\roman{Der}(A)$
where the notation $[\cdot,\cdot]$ and $\pi^{\sharp}$ is slightly abused.
For example, when $N$ is an ordinary smooth manifold
and $f$ a smooth function on $N$,
in local coordinates $(x_1,\dots,x_n)$, 
we have the formal differential $\alpha = df - \sum \partial_j f dx_j$,
and the Taylor theorem entails that this formal differential vanishes at every
point of $N$.
In \cite\poiscoho, we have pointed out 
that, in this particular case, the fact that the 
$(\Bobb R,A)$-Lie algebra structure on
$D_A$ descends to 
$\Omega^1(N)$
amounts to the nowadays familiar 
Lie algebroid structure on the cotangent bundle 
of a smooth Poisson manifold $N$.
{\sl When $N$ is a space with singularities,
the above description in terms
of the quotient
$\Omega^1(N)$
of $D_A$ by the differentials which vanish at every point of $N$ is
more general, though, and cannot be given in terms of Lie algebroids\/}.

\smallskip\noindent
(1.4) {\sl The stratified symplectic structure on the closure of any stratum\/.}
Let $(N,C^{\infty}(N),\{\cdot,\cdot\})$ be a stratified symplectic space.
The closure $\overline Y$ of any stratum $Y$
of $N$ inherits a stratified symplectic structure
$(C^{\infty}(\overline Y),\{\cdot,\cdot\}^{\overline Y})$
in the following fashion:
Let 
$C^{\infty}(\overline Y)$ be the algebra of continuous functions
on $\overline Y$ which arise from restriction 
to $\overline Y$
of functions in
$C^{\infty}(N)$.
Since 
the inclusion $Y \subseteq N$ is a Poisson map
the ideal of functions in
$C^{\infty}(N)$
which vanish on
$Y$ and hence on $\overline Y$
is a Poisson ideal.
Consequently
the Poisson structure on
$C^{\infty}(N)$
descends to a Poisson structure
on $C^{\infty}(\overline Y)$ which we denote by $\{\cdot,\cdot\}^{\overline Y}$;
thus 
$(Y,C^{\infty}(\overline Y),\{\cdot,\cdot\}^{\overline Y})$
is a stratified symplectic space,
and we refer to its structure as being {\it induced from the stratified
symplectic structure of \/} $N$.
The projection mapping
from $C^{\infty}(N)$ to $C^{\infty}(\overline Y)$
is plainly a morphism 
$\phi \colon C^{\infty}(N) \to C^{\infty}(\overline Y)$
of Poisson algebras which,
in view of the Addendum 3.8.4 in \cite\poiscoho,
induces a morphism 
$$
(C^{\infty}(N),D_{\{\cdot,\cdot\}}) @>>> 
(C^{\infty}(\overline Y),D_{\{\cdot,\cdot\}^{\overline Y}})
\tag1.4.1
$$
of Lie-Rinehart algebras.
For the record, we spell out the following, the proof of which is 
straightforward and left to the reader.

\proclaim{Proposition 1.4.2}
Given a stratified symplectic space $N$,
for any stratum $Y$, 
the closure $\overline Y$ being endowed with the induced
stratified symplectic 
structure explained above,
the morphism {\rm (1.4.1)}
passes to a morphism
$$
(\phi,\phi_*)\colon (C^{\infty}(N),\Omega^1(N)_{\{\cdot,\cdot\}}) @>>> 
(C^{\infty}(\overline Y),
\Omega^1(\overline Y)_{\{\cdot,\cdot\}^{\overline Y}})
\tag1.4.3
$$
of Lie-Rinehart algebras.
\endproclaim

\smallskip\noindent
(1.5) {\sl Stratified prequantum modules\/.}
In the constructions explained in (1.2) above,
we now replace
the 
$(\Bobb R,A)$-Lie algebra
$D_{\{\cdot,\cdot\}}$
with the
$(\Bobb R,A)$-Lie algebra
$\Omega^1(N)_{\{\cdot,\cdot\}}$ 
and, accordingly, we replace
the extension (1.2.1)
of $(\Bobb R,A)$-Lie algebras
with the corresponding extension
$$
0 
@>>> 
A
@>>> 
\overline L_{\{\cdot,\cdot\}}
@>>> 
\Omega^1(N)_{\{\cdot,\cdot\}}
@>>> 
0
\tag1.5.1
$$
of $(\Bobb R,A)$-Lie algebras.
Given an
$(A\otimes \Bobb C)$-module
$M$,
we refer to
an 
$(A,\overline L_{\{\cdot,\cdot\}})$-module
structure 
$$
\chi
\colon 
\overline L_{\{\cdot,\cdot\}}
@>>>
\roman{End}_{\Bobb R}(M)
\tag1.5.2
$$
on $M$ 
as a {\it geometric prequantum module\/} structure
or, more simply, as a
{\it prequantum module\/} structure, provided
the composite of $\chi$
with the canonical epimorphism
from
$\overline L^{\roman{a}}_{\{\cdot,\cdot\}}$ onto
$\overline L_{\{\cdot,\cdot\}}$
is an algebraic prequantum module structure.
A pair $(M,\chi)$
consisting of 
an $(A\otimes \Bobb C)$-module $M$ and a (geometric) prequantum module structure
will be referred to as a ({\it geometric\/}) {\it prequantum module\/}
(for the stratified symplectic space $(N,C^{\infty}(N),\{\cdot,\cdot\})$.
\smallskip
In particular, let $(N,\sigma)$ be a quantizable
symplectic manifold.
Consider
the 2-form $\pi$
induced by the symplectic Poisson structure and,
as in (1.2) above, write
the adjoint of $\pi$
as
$\pi^{\sharp}\colon \Omega^1(N) \to \roman{Vect}(N)$.
Under the circumstances of (1.1) above,
let $M= \Gamma^{\infty}(\zeta)$;
the assignments
$$
\chi_{\nabla}(a,0) = i\,a\,\roman{Id}_M,
\quad
\chi_{\nabla}(0,\alpha) = \nabla_{\pi^{\sharp}(\alpha)},
\quad a \in A,\ \alpha \in \Omega^1(N),
$$
yield 
a geometric prequantum module structure 
$$
\chi_{\nabla} 
\colon
\overline L_{\{\cdot,\cdot\}} @>>>
\roman{End}_{\Bobb C}(M)
\subseteq 
\roman{End}_{\Bobb R}(M)
\tag1.5.3
$$
for
$(A,\{\cdot,\cdot\})$.
This is just
the ordinary prequantization
construction
in another guise.
In fact,
under the adjoint $\pi^{\sharp}$ 
from 
$\Omega^1(N)$ to
$\roman{Vect}(N)$,
the 
prequantization
formula (1.2.8) passes to the more usual prequantization
formula (1.1.1).
Occasionally we shall refer to a prequantum module structure
of the kind 
(1.5.3)
as {\it smooth\/}.
\smallskip
Let $N$ be a stratified symplectic space,
with stratified symplectic Poisson algebra
$(C^{\infty}(N),\{\cdot,\cdot\})$.
For each stratum $Y$, let
$(C^{\infty}(Y),\{\cdot,\cdot\}^Y)$
be its ordinary {\it smooth\/} symplectic Poisson structure,
and let
$$
0
@>>>
C^{\infty}(Y)
@>>>
\overline L_{\{\cdot,\cdot\}^Y}
@>>>
\Omega^1(Y)_{\{\cdot,\cdot\}^Y}
@>>>
0
\tag1.5.4
$$
be the corresponding extension (1.5.1) of $(\Bobb R,A)$-Lie algebras
where 
$\Omega^1(Y)$
is the (projective)
$C^{\infty}(Y)$-module of ordinary 1-forms on $Y$.
We define a {\it stratified prequantum module\/}
for $N$ to consist of
\newline\noindent
--- a (geometric) prequantum module $(M,\chi)$ for
$(C^{\infty}(N),\{\cdot,\cdot\})$
having the property that, for any stratum $Y$,
the canonical
linear 
map 
of complex vector spaces
from $M_{\overline Y} =  C^{\infty}(\overline Y) \otimes_{C^{\infty}(N)} M$ to 
$M_Y =  C^{\infty}(Y) \otimes_{C^{\infty}(N)} M$ is injective,
together with,
\newline\noindent
--- for each stratum $Y$, a prequantum module structure
$\chi_Y$ for
$(C^{\infty}(Y),\{\cdot,\cdot\}^Y)$
on the induced module
$M_Y =  C^{\infty}(Y) \otimes_{C^{\infty}(N)} M$
in such a way that the canonical
linear 
map 
of complex vector spaces
from $M$ to $M_Y$
is a
morphism 
of prequantum modules
from $(M,\chi)$ to $(M_Y,\chi_Y)$.
\smallskip
Here \lq morphism 
of prequantum modules
from $(M,\chi)$ to $(M_Y,\chi_Y)$\rq\ 
means that
\linebreak
(i) the canonical
linear 
map 
of complex vector spaces
from $M$ to $M_Y$, (ii) the adjoints 
$\chi^{\sharp}$
and
$\chi_Y^{\sharp}$
of the structure maps,
and (iii) the morphism
from
$\overline L_{\{\cdot,\cdot\}}$
to
$\overline L_{\{\cdot,\cdot\}^Y}$
induced by the restriction map,
make commutative the diagram
$$
\CD
\phantom{\,}\overline L_{\{\cdot,\cdot\}}\otimes M \phantom{\{\,\,\,\cdot\}}
@>{\chi^{\sharp}}>>
M
\\
@VVV
@VVV
\\
\overline L_{\{\cdot,\cdot\}^Y} \otimes M_Y \phantom{\{\}\,}
@>{\chi_Y^{\sharp}}>>
\phantom{\,}M_Y .
\endCD
\tag1.5.5
$$
Occasionally we shall denote a stratified prequantum
module structure by $(\chi, \{\chi_Y\})$
where $Y$ runs through the strata, or 
sometimes more simply just by $\chi$,
with an abuse of notation.
For intelligibility, we note that,
on a stratum $Y$,
the induced module
$M_Y =  C^{\infty}(Y) \otimes_{C^{\infty}(N)} M$
will often come down to the space of sections of
an ordinary smooth complex (perhaps V-) line bundle
and, perhaps with a grain of salt, a prequantum module structure 
$\chi_Y$ for
$(C^{\infty}(Y),\{\cdot,\cdot\}^Y)$
on  $M_Y$
will then come down to ordinary prequantization,
cf. (1.5.3).

\proclaim{Theorem 1.5.6}
Let $(N,C^{\infty}(N),\{\cdot,\cdot\})$
be a stratified symplectic space,
let $Y$ be a stratum of $N$, 
let $(\overline Y,C^{\infty}(\overline Y),\{\cdot,\cdot\}^{\overline Y})$
be the induced stratified symplectic structure on
the closure $\overline Y$ of $Y$,
cf. {\rm (1.4)},
and let $(M,\chi)$ be a stratified prequantum module for
$(N,C^{\infty}(N),\{\cdot,\cdot\})$.
The induced $C^{\infty}(\overline Y)$-module
$M_{\overline Y} =  C^{\infty}(\overline Y) \otimes_{C^{\infty}(N)} M$
inherits a
stratified prequantum module structure
$$
\chi_{\overline Y}\colon
\overline L_{\{\cdot,\cdot\}^{\overline Y}}
@>>>
\roman{End}_{\Bobb R}(M_{\overline Y})
\tag1.5.7
$$
for 
$(\overline Y,C^{\infty}(\overline Y),\{\cdot,\cdot\}^{\overline Y})$
in such a way that the canonical
linear 
map 
of complex vector spaces
from $M$ to $M_{\overline Y}$
is a
morphism 
of stratified prequantum modules
from $(M,\chi)$ to $(M_{\overline Y},\chi_{\overline Y})$.
\endproclaim

\demo{Proof}
The morphism (1.4.3)
of Lie-Rinehart algebras plainly induces a morphism
$$
(\phi,\phi_{\flat})\colon 
(C^{\infty}(N),\overline L_{\{\cdot,\cdot\}}) 
@>>> 
(C^{\infty}(\overline Y),\overline L_{\{\cdot,\cdot\}^{\overline Y}})
\tag1.5.8
$$
of Lie-Rinehart algebras in such a way that
the restriction
$$
(C^{\infty}(\overline Y),\overline L_{\{\cdot,\cdot\}^{\overline Y}})
@>>>
(C^{\infty}(Y),
\overline L_{\{\cdot,\cdot\}^{Y}})
\tag1.5.9
$$
to the stratum $Y$, combined with (1.5.8),
amounts to the restriction morphism
$$
(\phi,\phi_{\flat})\colon (C^{\infty}(N),
\overline L_{\{\cdot,\cdot\}}) @>>> 
(C^{\infty}(Y),
\overline L_{\{\cdot,\cdot\}^{Y}})
\tag1.5.10
$$
from $N$ to the stratum $Y$.
To isolate what we must precisely prove,
write
$(A,L) = (C^{\infty}(N),\overline L_{\{\cdot,\cdot\}})$, 
$(A',L')=(C^{\infty}(\overline Y),\overline L_{\{\cdot,\cdot\}^{\overline Y}})$,
$M_{\overline Y} =M' = A' \otimes _A M$, and write the Poisson structures
on $A$ and $A'$ as $\{\cdot,\cdot\}$ and $\{\cdot,\cdot\}'$, respectively.
Let $I$ be the ideal of functions
in $C^{\infty}(N)$, necessarily a Poisson ideal (as we have already observed),
which vanish on $Y$ and hence on $\overline Y$.
By construction, $A'$ is canonically isomorphic to  $A/I$.
We must prove that the $L$-action 
$\chi \colon L \to
\roman{End}_{\Bobb R}(M)$
on $M$  
passes to an $L'$-action
$\chi' \colon L' \to\roman{End}_{\Bobb R}(M')$ 
on $M'$.
In view of standard arguments
from the theory of formal differentials, the canonical surjection
from $D_A$ to $D_{A'}$ fits into the exact sequence
$$
I/I^2 @>>> A' \otimes _A D_A @>>> D_{A'} @>>> 0
\tag1.5.11
$$
of $A'$-modules,
in fact, $(\Bobb R,A')$-Lie algebras,
where the unlabelled left-hand arrow is given by
the assignment to $f \in I$ of $1 \otimes df$.
This exact sequence, in turn, lifts to an exact sequence
$$
I/I^2 @>>> A' \otimes _A \overline L^{\roman{a}}_{\{\cdot,\cdot\}}
@>>> \overline L^{\roman{a}}_{\{\cdot,\cdot\}'} @>>> 0
$$
of $A'$-modules,
even $(\Bobb R,A')$-Lie algebras,
and dividing out the differentials which vanish at every point
(of $\overline Y$) we obtain an exact sequence
$$
I/I^2 @>>> A' \otimes _A L
@>>> L' @>>> 0
$$
of $(\Bobb R,A')$-Lie algebras.
Since $M' = A' \otimes _A M$,
the $(A,L)$-module structure on $M$ passes to an
$(A',A' \otimes _A L)$-module structure on $M'$.
Thus we must prove that, whenever $f \in I$,
that is, whenever $f$ is a function in 
$C^{\infty}(N)$ which vanishes on $Y$,
the element $1 \otimes df$ of $A' \otimes _A L$
acts trivially on $M'$.
This may be seen as follows:
\smallskip
Since $(M,\chi)$ is a stratified prequantum module,
in view of the definition,
with reference to the stratum $Y$ 
(which is an ordinary smooth manifold),
the induced module
$M_Y =  C^{\infty}(Y) \otimes_{C^{\infty}(N)} M$
has a prequantum module structure
$\chi_Y$ for
$(C^{\infty}(Y),\{\cdot,\cdot\}^Y)$
in such a way that the canonical
linear 
map 
of complex vector spaces
from $M$ to $M_Y$
is a
morphism 
of prequantum modules
from $(M,\chi)$ to $(M_Y,\chi_Y)$.
The $C^{\infty}(Y)$-module underlying the latter may be written as
$$
M_Y =  C^{\infty}(Y) \otimes_{C^{\infty}(N)} M
\cong
C^{\infty}(Y) \otimes_{A'} M',
$$ 
and the induced morphism from
$M' \cong M_{\overline Y}$ to $M_Y$ is compatible with the actions,
with respect to the induced morphism
of Lie-Rinehart algebras
from $(A',A' \otimes _A L)$
to
$(C^{\infty}(Y), \Omega^1(Y)_{\{\cdot,\cdot\}^Y})$;
that is to say: these morphisms make the diagram
$$
\CD
\phantom{\,}(A' \otimes _A L)\otimes M' \phantom{\{\,\,\,\cdot\}}
@>>>
M'
\\
@VVV
@VVV
\\
\overline L_{\{\cdot,\cdot\}^Y} \otimes M_Y \phantom{\{\}\,}
@>>>
\phantom{\,}M_Y 
\endCD
$$
commutative.
Consequently the element $1 \otimes df$ of $A' \otimes _A L$
acts trivially on 
$M_Y$ since it acts thereupon through the map from
$A' \otimes _A L$
to 
$\overline L_{\{\cdot,\cdot\}^Y}$
where it becomes trivial,
and thence $1 \otimes df$
acts trivially on $M' \cong M_{\overline Y}$
since, by definition, the canonical map of complex vector spaces from
$M_{\overline Y}$ to $M_Y$ is required to be injective. \qed
\enddemo

\noindent
(1.6) {\sl Costratified prequantum spaces\/.}
Let $(N,C^{\infty}(N),\{\cdot,\cdot\})$
be a stratified symplectic space.
Under the circumstances of (1.5.6),
when $Y$ runs through the strata of $N$, 
we will refer to the system 
$$
\left(M_{\overline Y},
\chi_{\overline Y}\colon \overline L_{\{\cdot,\cdot\}^{\overline Y}}
@>>> \roman{End}_{\Bobb R}(M_{\overline Y})\right)
\tag1.6.1
$$
of prequantum modules, together with,
for every
pair of strata $Y,Y'$ such that
$Y' \subseteq \overline Y$,
the induced morphism
$$
\left(M_{\overline Y},\chi_{\overline Y}\right)
@>>>
\left(M_{\overline Y'},\chi_{\overline Y'}\right)
\tag1.6.2
$$
of prequantum modules,
as a {\it costratified prequantum space\/}.
More formally: Consider the category $\Cal C_N$ whose objects are the strata
of $N$ and whose morphisms are the inclusions $Y' \subseteq \overline Y$.
We define a 
{\it costratified complex vector space\/} on $N$ to be a contravariant functor
from $\Cal C_N$ to the category of complex vector spaces,
and a
{\it costratified prequantum space\/} on $N$ to be a 
costratified complex vector space
together with a compatible system of prequantum module structures.
The discussion in (1.4) and (1.5) may be summarized by saying that
a stratified prequantum module  $(M,\chi)$ 
for 
$(N,C^{\infty}(N),\{\cdot,\cdot\})$
determines a costratified prequantum space on $N$:
For every stratum $Y$,
let $(\overline Y,C^{\infty}(\overline Y),\{\cdot,\cdot\}^{\overline Y})$
be the induced stratified symplectic structure on
the closure $\overline Y$ of $Y$ given in (1.4),
and let 
$(C^{\infty}(\overline Y),\overline L_{\{\cdot,\cdot\}^{\overline Y}})$
be the corresponding Lie-Rinehart algebra.
Then the 
assignment to a stratum $Y$ of
the induced stratified prequantum module
$\left(M_{\overline Y},\chi_{\overline Y}\right)$
is a functor on
$\Cal C_N$
in an obvious fashion;
in particular,
whenever $Y$ and $Y'$ are two strata such that $Y' \subseteq \overline Y$,
restriction yields a morphism
$$
\left(M_{\overline Y},\chi_{\overline Y}\right)
@>>> \left(M_{\overline {Y'}},\chi_{\overline {Y'}}\right)
$$
of stratified prequantum modules.
The system which encompasses {\it all\/}
these restriction morphisms is the corresponding
costratified prequantum space.

\medskip\noindent {\bf 2. Prequantum modules and reduction}
\smallskip\noindent
Consider an ordinary smooth symplectic manifold $N$,
acted upon by a compact Lie group $G$ in a hamiltonian fashion
with momentum mapping $\mu \colon N \to \fra g^*$
where $\fra g$ is the Lie algebra of $G$.
Suppose in addition that the symplectic manifold $N$ is quantizable,
let $\zeta\colon \Lambda \to N$ be a prequantum bundle,
and suppose that
the hamiltonian $G$-action on $N$
lifts to an action on $\zeta$ preserving
the connection.
When $G$ is connected,
this additional assumption is (well known to be) redundant
(since the action is hamiltonian) and
it will suffice to
replace $G$ with a suitable covering group 
if need be so that 
the hamiltonian action on $N$
lifts to an action of the covering group on $\zeta$ preserving
the connection,
cf. e.~g. \cite\kostaone; we then
work with this covering group
which we continue to denote by $G$.
Thus, 
with this preparation out of the way,
$G$ acts on $\zeta$.
Let $M= \Gamma^{\infty}(\zeta)$
(the space of smooth complex valued sections of $\zeta$), endowed with the
smooth prequantum module structure (1.5.3), which we now write
as $\chi \colon \overline L_{\{\cdot,\cdot\}} \to 
\roman{End}_{\Bobb R}(M)$.
In this paper we will only explore the case of zero value of the momentum
mapping. The case of non-zero value would require much more effort and,
in particular, the connection preserving
prequantum lift of the action might be inappropriate.
\smallskip
The reduced space $N^{\roman{red}} = \mu^{-1}(0) \big / G$
is well known to be stratified by orbit types.
As usual,
the superscript 
\lq\lq\ $-^G$\ \rq\rq\ will refer to $G$-invariants.
Let $I$ be the ideal of smooth functions on $N$
which vanish on the zero locus
$\mu^{-1}(0)$, and let
$C^{\infty}(N^{\roman{red}})=(C^{\infty}(N))^G /I^G$,
the algebra of smooth $G$-invariant functions, divided out by the ideal
of smooth
$G$-invariant functions
which vanish on the zero locus.
This is an algebra of continuous functions on
$N^{\roman{red}}$ in an obvious fashion such that, 
restricted to any stratum, these functions are smooth;
in other words,
$C^{\infty}(N^{\roman{red}})$
is a smooth structure on
$N^{\roman{red}}$.
As observed in {\smc Arms-Cushman-Gotay} \cite\armcusgo, 
cf. the proof of Theorem 1 in \cite\armcusgo, 
the Noether theorem entails that 
the ideal $I^G$ of $G$-invariant functions
which vanish on
$\mu^{-1}(0)$ is a Poisson ideal in
the Poisson algebra $(C^{\infty}(N))^G$
of smooth $G$-invariant functions, that is,
the symplectic Poisson structure on 
$C^{\infty}(N)$ descends to a Poisson structure 
$\{\cdot,\cdot\}^{\roman{red}}$ on
$C^{\infty}(N^{\roman{red}})$.
\smallskip
Consider the 
sheaf or complex V-line bundle
(also called orbi bundle)
$$
\zeta^{\roman{red}}
\colon 
\Lambda^{\roman{red}} = \left(\Lambda\big | \mu^{-1}(0)\right)\big /G
@>>>
N^{\roman{red}}.
$$
As a sheaf (in the category of spaces and ordinary
continuous functions)
it may be written as the direct image 
$q_*\left(\zeta \big | \mu^{-1}(0)\right)$
where $q \colon \mu^{-1}(0) \to N^{\roman{red}}$
refers to the projection.
On each stratum,
cf. e.~g. \cite{\guistetw,\,\sjamatwo},
this V-line bundle
restricts to an ordinary smooth complex
V-line bundle.
The sheaf (or V-line bundle)
$\zeta$
gives rise to a prequantization construction
for the reduced Poisson algebra, in the following fashion:
We will write
the space
of sections
of $\zeta$ that vanish on the zero locus
of $\mu$
as $IM$.
We mention in passing that
a partition of unity argument shows that
any section
of $\zeta$ that vanishes on the zero locus
may be written as a sum
$\sum \beta_j \zeta_j$
where 
the $\zeta_j$'s are just sections of $\zeta$
and where the $\beta_j$'s are functions in $I$, that is, functions
which vanish on the zero locus $\mu^{-1}(0)$.
This justifies the notation.
We shall not exploit this observation, though.
Let $M^{\roman{red}}= M^G\big /(IM)^G$, that is, the space $M^G$ of 
ordinary 
smooth $G$-invariant sections of $\zeta$, modulo the subspace
$(IM)^G$ of 
smooth $G$-invariant sections
that vanish on the zero locus $\mu^{-1}(0)$.
By construction, $M^{\roman{red}}$ 
may canonically be identified with
a space of continuous sections
of $\zeta^{\roman{red}}$ which, on each stratum, restrict to a smooth section,
and $M^{\roman{red}}$
inherits a $C^{\infty}(N^{\roman{red}},\Bobb C)$-module structure.
Accordingly,
we will occasionally denote $M^{\roman{red}}$
by $\Gamma^{\infty}(\zeta^{\roman{red}})$
when $M^{\roman{red}}$
is viewed being endowed with this
$C^{\infty}(N^{\roman{red}},\Bobb C)$-module structure.
As a side remark we note that,
since associating to
an ordinary complex vector bundle 
on $N^{\roman{red}}$
its space of continuous sections is an equivalence
of categories from
complex vector bundles to projective
modules over the algebra $C(N^{\roman{red}},\Bobb C)$
of continuous functions on
$N^{\roman{red}}$,
when $M^{\roman{red}}$ is a projective
$C^{\infty}(N^{\roman{red}},\Bobb C)$-module
the reduced V-line bundle $\zeta^{\roman{red}}$
is an ordinary line bundle.
The converse is presumably true as well;
it would rely on the corresponding equivalence of categories
spelled out for
$C^{\infty}(N^{\roman{red}},\Bobb C)$-modules
but details have not been worked out yet.

\proclaim{Theorem 2.1}
Given an ordinary $G$-equivariant prequantum bundle
$\zeta$ on the smooth hamiltonian $G$-space $(N,\mu)$,
the prequantum module structure
$\chi$ of 
$M=\Gamma^{\infty}(\zeta)$ determines a stratified prequantum module
structure 
$\chi^{\roman{red}}$
on $M^{\roman{red}}$
{\rm (made explicit in (2.9) below)}
for the stratified symplectic space
$(N^{\roman{red}},C^{\infty}(N^{\roman{red}}),\{\cdot,\cdot\}^{\roman{red}})$.
\endproclaim

\demo{Proof}
Let $A=C^{\infty}(N)$, 
write $\{\cdot,\cdot\}^G$ for the induced Poisson structure on $A^G$,
and let $\nabla$ be the operator of covariant derivative
determined by the connection on the prequantum bundle $\zeta$.
The corresponding extension (1.2.1)
of $(\Bobb R,A^{G})$-Lie algebras 
may be written as
$$
0 
@>>> 
A^{G}
@>>> 
\overline L^{\roman{a}}_{\{\cdot,\cdot\}^{G}}
@>>> 
D_{\{\cdot,\cdot\}^{G}}
@>>> 
0
\tag2.2
$$
and, by naturality,
the inclusion
 $(A^{G},\{\cdot,\cdot\}^{G}) \to 
(A,\{\cdot,\cdot\})$ 
of Poisson algebras
induces a commutative diagram
$$
\CD
0 
@>>> 
A^{G}
@>>> 
\overline L^{\roman{a}}_{\{\cdot,\cdot\}^{G}}
@>>> 
D_{\{\cdot,\cdot\}^{G}}
@>>> 
0\phantom{.}
\\
@.
@VVV
@VVV
@VVV
@.
\\
0 
@>>> 
A
@>>> 
\overline L^{\roman{a}}_{\{\cdot,\cdot\}}
@>>> 
D_{\{\cdot,\cdot\}}
@>>> 
0
\endCD
\tag2.3
$$
of extensions of Lie-Rinehart algebras;
notice that here the algebra variable changes 
whence we cannot refer merely to $(\Bobb R,A)$-Lie algebras.
The composite of $\chi$ 
(more precisely, of the corresponding algebraic prequantum module structure)
with
the induced map
from
$\overline L^{\roman{a}}_{\{\cdot,\cdot\}^{G}}$
to 
$\overline L^{\roman{a}}_{\{\cdot,\cdot\}}$
yields an
$(A^{G},\overline L^{\roman{a}}_{\{\cdot,\cdot\}^{G}})$-module
structure on
$M$.
\smallskip
By symmetry,
the $\overline L^{\roman{a}}_{\{\cdot,\cdot\}^{G}}$-action
on $M$ preserves $M^G$
since, given a $G$-equivariant section $\eta$ 
of $\zeta$ and a $G$-equivariant vector field $X$ on $N$,
the section $\nabla_X \eta$ of $\zeta$ is also $G$-equivariant.
Hence
$M^G$ is a submodule
of $M$,
for the $(A^{G},\overline L^{\roman{a}}_{\{\cdot,\cdot\}^{G}})$-module
structure.
We write
$\chi^G \colon \overline L^{\roman{a}}_{\{\cdot,\cdot\}^{G}} 
\to \roman{End}_{\Bobb R}(M^G)$
for the 
resulting
algebraic prequantum module 
structure 
for
$(A^{G},\{\cdot,\cdot\}^{G})$
on $M^G$.
We now assert that
this 
$(A^{G},\overline L^{\roman{a}}_{\{\cdot,\cdot\}^{G}})$-module
structure $\chi^{G}$ 
preserves 
$(IM)^G \subseteq M^{G}$.
Since we already know that $\chi^{G}$  preserves
$M^{G}$ it will suffice to see that
it preserves $IM$.
\smallskip
In order to justify the last claim, 
let $\Lambda^\times \subseteq \Lambda$ be the subspace of non-zero vectors,
so that
$\zeta$, restricted to
$\Lambda^\times$,
is the corresponding principal $\Bobb C^\times$-bundle
which, abusing the notation $\zeta$, we 
denote by $\zeta\colon \Lambda^\times \to N$ as well.
The identity
$$
h(\zeta(z)) =h^\sharp (z) z,\quad z \in \Lambda^\times,
$$
is well known to establish
an isomorphism between
the space $M$ of sections $h$ of $\zeta$ and
the space of functions $h^\sharp \colon 
\Lambda^\times \to \Bobb C$
which are {\it equivariant\/} in the sense that
$$
h^\sharp (cz) = c^{-1} h^\sharp(z),\quad z \in \Lambda^\times,\ c
\in \Bobb C^\times.
$$
For any vector field $X$ on $N$,
the covariant derivative of a section $h$ in the direction of $X$
is the section
$\nabla_X h$ of $\zeta$ given by
$$
\nabla_X h(\zeta(z)) = (X^{\sharp} h^{\sharp})(z) z,\quad 
z \in \Lambda^\times,
$$
where 
$X^{\sharp}$ is the horizontal lift of $X$ 
(with reference to the connection corresponding to
the prequantum bundle structure)
to a vector field on
$\Lambda^\times$;
see e.~g. \cite\sniabook\ for details.
Let $f \in (C^{\infty}(N))^G$.
By the Noether theorem,
the momentum mapping $\mu \colon N \to \fra g^*$
is constant along the flow lines of the hamiltonian vector field
$X_f$ ($= \{\cdot,f\}$) of $f$
as well as 
along the flow lines of the vector field
$-X_f= \{f,\cdot\}$ which comes into play in the prequantization formula
(1.1.1).
Since each flow line
of
$X_f^{\sharp}$ in $\Lambda^\times$
is the (unique) horizontal lift
of a flow line 
of $X_f$ in $N$,
the composite
$$
\mu^{\sharp}= \mu \circ \zeta \colon  \Lambda^\times @>>> \fra g^*
$$
is constant along
the flow lines  of
$X_f^{\sharp}$ in $\Lambda^\times$.
\smallskip
Let $h^\sharp \colon 
\Lambda^\times \to \Bobb C$
be an equivariant function which vanishes on
$(\mu^{\sharp})^{-1}(0)$.
Thus 
$h^\sharp$
corresponds to a section $h$ of $\zeta$ 
which vanishes on $\mu^{-1}(0)$, i.~e.
$h \in IM$.
Let $q \in (\mu^{\sharp})^{-1}(0)$,
and
let $f \in (C^{\infty}(N))^G$.
Consider a flow line
$$
J @>>> \Lambda^\times,\quad t \mapsto \roman{exp}(tX^{\sharp}_f)q,
$$
of $X_f^{\sharp}$
having $q$ as its starting point
where $J$ is a suitable open interval containing $0$.
Since
$\mu^{\sharp}$
is constant along
the flow lines of $X_f^{\sharp}$,
$\mu^{\sharp}(
\roman{exp}(tX^{\sharp}_f)q) = 0$ for every $t \in J$.
Consequently
$$
(X_f^{\sharp} (h^\sharp))(q) = 
\frac d{dt}\left(h^\sharp(
\roman{exp}(tX^{\sharp}_f)q)\right)\big|_{t=0} =0, 
$$
that is, the vanishing of the section $h$
on $\mu^{-1}(0)$ implies that
of $\nabla_{X_f} h$ as well.
This observation entails that the 
$(A^{G},\overline L^{\roman{a}}_{\{\cdot,\cdot\}^{G}})$-module
structure $\chi^{G}$ 
preserves $IM \subseteq M$ since $f$ is an arbitrary $G$-invariant function.
For 
on the Lie-Rinehart generators
$(0,\alpha)$ 
of $\overline L_{\{\cdot,\cdot\}}= A \oplus \Omega^1(N)$
where $\alpha \in \Omega^1(N)$,
cf. (1.2.2) for the corresponding \lq\lq algebraic\rq\rq\ 
construction,
the prequantum module structure (1.5.3)
written there as $\chi_{\nabla}$,
is given by
$\chi_{\nabla}(0,\alpha) = \nabla_{\pi^{\sharp}(\alpha)}$,
and it suffices to take here as generators differentials
$\alpha$ of the kind $\alpha = dh$
where $h$ runs through (smooth) functions on $N$.
Likewise,
as an $(\Bobb R, A^G)$-Lie algebra,
$\overline L^{\roman{a}}_{\{\cdot,\cdot\}^{G}}$
is generated by differentials 
$\alpha = df$
of $G$-invariant functions $f$ on $N$,
more precisely, by the elements
$$
(0,df) \in 
A^G \oplus D_{\{\cdot,\cdot\}^G}
=\overline L^{\roman{a}}_{\{\cdot,\cdot\}^G}.  
$$
Consequently
the 
$(A^{G},\overline L^{\roman{a}}_{\{\cdot,\cdot\}^{G}})$-module
structure $\chi^{G}$ 
preserves $IM \subseteq M$
as asserted.
Since it also
preserves $M^G \subseteq M$,
it preserves $(IM)^G \subseteq M^G$ 
and therefore
induces an
$(A^{G},\overline L^{\roman{a}}_{\{\cdot,\cdot\}^{G}})$-module structure
$$
\overline L^{\roman{a}}_{\{\cdot,\cdot\}^{G}}
\longrightarrow
\roman{End}_{\Bobb R}(M^{\roman{red}})
\tag2.4
$$
on
$M^{\roman{red}} =M^{G}/(IM)^G$.
\smallskip
The obvious epimorphism from
$D_{\{\cdot,\cdot\}^{G}}$
onto
$D_{\{\cdot,\cdot\}^{\roman{red}}}$
and that from
$\overline L^{\roman{a}}_{\{\cdot,\cdot\}^{G}}$
onto $\overline L^{\roman{a}}_{\{\cdot,\cdot\}^{\roman{red}}}$
together give rise to  {\it induced\/}
$(\Bobb R,A^{\roman{red}})$-Lie algebra structures   on
$A^{\roman{red}}\otimes _{A^{G}} D_{\{\cdot,\cdot\}^{G}}$
and
$A^{\roman{red}}\otimes _{A^{G}} \overline L^{\roman{a}}_{\{\cdot,\cdot\}^{G}}$,
respectively;
see Section 1 in \cite\poiscoho\ for 
details on induced structures.
Application of the functor
$A^{\roman{red}}\otimes _{A^{G}} - $
to (2.2) then yields the extension
$$
0 
@>>> 
A^{\roman{red}}
@>>> 
A^{\roman{red}}\otimes _{A^{G}}\overline L^{\roman{a}}_{\{\cdot,\cdot\}^{G}}
@>>> 
A^{\roman{red}}\otimes _{A^{G}} D_{\{\cdot,\cdot\}^{G}}
@>>> 
0
\tag2.5
$$
of
$(\Bobb R,A^{\roman{red}})$-Lie algebras.
It is clear that
(2.4) factors through
the obvious epimorphism from
$\overline L^{\roman{a}}_{\{\cdot,\cdot\}^{G}}$
onto 
$A^{\roman{red}}\otimes _{A^{G}}\overline L^{\roman{a}}_{\{\cdot,\cdot\}^{G}}$.
\smallskip
We now assert that
(2.4) factors even through
the epimorphism from $\overline L^{\roman{a}}_{\{\cdot,\cdot\}^{G}}$
onto $\overline L^{\roman{a}}_{\{\cdot,\cdot\}^{\roman{red}}}$
induced by the epimorphism of Poisson algebras
from
$(A^{G},\{\cdot,\cdot\}^{G})$
to
$(A^{\roman{red}},\{\cdot,\cdot\}^{\roman{red}})$.
Indeed,
by naturality,
the 
inclusion
of Poisson algebras
from
$(A^{G},\{\cdot,\cdot\}^{G})$ into $(A,\{\cdot,\cdot\})$ 
induces the commutative diagram
$$
\CD
0 
@>>> 
A^{G}
@>>> 
\overline L^{\roman{a}}_{\{\cdot,\cdot\}^{G}}
@>>> 
D_{\{\cdot,\cdot\}^{G}}
@>>> 
0\phantom{.}
\\
@.
@VVV
@VVV
@VVV
@.
\\
0 
@>>> 
A^{\roman{red}}
@>>> 
A^{\roman{red}}\otimes _{A^{G}}\overline L^{\roman{a}}_{\{\cdot,\cdot\}^{G}}
@>>> 
A^{\roman{red}}\otimes _{A^{G}} D_{\{\cdot,\cdot\}^{G}}
@>>> 
0\phantom{.}
\\
@.
@VVV
@VVV
@VVV
@.
\\
0 
@>>> 
A^{\roman{red}}
@>>> 
\overline L^{\roman{a}}_{\{\cdot,\cdot\}^{\roman{red}}}
@>>> 
D_{\{\cdot,\cdot\}^{\roman{red}}}
@>>> 
0.
\endCD
\tag2.6
$$
By the general theory 
of K\"ahler differentials,
the obvious morphism 
of $ A^{\roman{red}}$-modules
from
$A^{\roman{red}} \otimes _{A^{G}}  D_{A^{G}}$
to
$D_{A^{\roman{red}}}$
fits into
an exact sequence
$$
 I^G/ (I^G)^2
\longrightarrow
 A^{\roman{red}} \otimes _{A^{G}}  D_{A^{G}}
\longrightarrow
 D_{A^{\roman{red}}}
\longrightarrow
0
\tag2.7
$$
of $A^{\roman{red}}$-modules
where the first arrow 
is given by the association
$$
a\  \roman{mod}\  (I^G)^2
\longmapsto
1 \otimes da \in  A^{\roman{red}} \otimes _{A^{G}}  D_{A^{G}},
$$
cf. (1.5.11) above.
Now we claim that,
for $f \in I^G$ and $h \in M^{G}$, we have
$$
df(h) =\nabla_{\pi^{\sharp}(df)}(h)\in (IM)^G
$$
where 
$df$ is viewed as the element 
$(0,df)$
of
$\overline L^{\roman{a}}_{\{\cdot,\cdot\}^{G}} 
= A^{G} \oplus D_{\{\cdot,\cdot\}^{G}}$.
To justify this claim we recall first that, as already noted above,
for any $q \in \mu^{-1}(0)$,
$$
\mu (\roman{exp}(tX_f)q) = \mu(q) = 0
$$
whence the flow line
of $X_f$ starting at $q$ lies in
$\mu^{-1}(0)$, in fact, in a stratum thereof.
Given a smooth $G$-invariant function $k$ on $N$,
by the Noether theorem,
the momentum mapping $\mu$ is constant along the trajectories
of $X_k$. Hence,
given $q \in \mu^{-1}(0)$,
since $\mu(\roman{exp}(tX_k)q) = \mu(q)$,
we have
$$
\{f,k\}(q) =-(X_f(k)) (q)=  
(X_k(f))(q) = \frac d{dt}\left(f(\roman{exp}(tX_k)q)\right) \big|_{t=0} = 0
$$
Hence
$$
(X_f(k))(q) =  
\frac d{dt}\left(k(\roman{exp}(tX_f)q)\right) \big|_{t=0} = 0 .
$$
Our aim is to prove
that, likewise, for
$q^{\sharp} \in \zeta^{-1}(q)$,
$$
(X_f^{\sharp} (h^\sharp))(q^{\sharp}) 
=
\frac d{dt}\left(h^\sharp(\roman{exp}(tX^\sharp_f)q^\sharp)\right)
\big|_{t=0} = 0
$$
where $h^\sharp$ is the smooth $G$-invariant function on $\Lambda^\times$ 
which corresponds to the smooth $G$-equivariant
section $h$ of $\zeta$ as explained above.
To this end, we refine the reasoning in the following fashion:
Let $f$ be any smooth function on $N$ which is zero on $\mu^{-1}(0)$,
at this stage not necessarily $G$-invariant and,
as before, let $q \in \mu^{-1}(0)$.
Then, for any $Y \in \roman{ker}(d\mu_q\colon \roman T_qN \to \fra g^*)$,
we have
$$
\omega(X_f\big|_q, Y) = Y f = (Y,df(q)) = 0.
$$
Since the annihilator of
$\roman{ker}(d\mu_q)$ coincides with
the tangent space
$\roman T_q Gq \subseteq \roman T_qN$ of the $G$-orbit $Gq \subseteq N$,
we conclude that
$X_f\big|_q \in \roman T_q Gq$.
Hence, given any smooth $G$-invariant function $k$ on $N$,
$$
(X_f(k)) (q) = X_f\big|_q (k) = 0.
$$
Since the connection on $\zeta$ is $G$-equivariant,
the horizontal lift
$X^{\sharp}_f\big|_{q^{\sharp}}$ 
of $X_f\big|_q$
lies in  
$\roman T_{q^{\sharp}} Gq^{\sharp} \subseteq
\roman T_{q^{\sharp}}  \Lambda^\times$. 
Consequently
$$
(X_f^{\sharp} (h^\sharp))(q^{\sharp}) 
=
X_f^{\sharp}\big|_{q^{\sharp}} (h^\sharp) = 0,
$$
for $h^\sharp$ is $G$-invariant.
Since 
$q \in \mu^{-1}(0)$ is arbitrary, we conclude that
$$
\nabla_{X_f} h
=
\nabla_{\pi^{\sharp}(df)}(h)
$$
vanishes on $\mu^{-1}(0)$, that is, that
$\nabla_{\pi^{\sharp}(df)}$ lies in $IM$.
By symmetry,
if, in addition, $f$ is $G$-invariant,
$\nabla_{\pi^{\sharp}(df)}$ 
is $G$-equivariant whence
$\nabla_{\pi^{\sharp}(df)}(h)$ lies in $(IM)^G$
as asserted.
Consequently
(2.4)
factors through the epimorphism from
$\overline L^{\roman{a}}_{\{\cdot,\cdot\}^{G}}$
onto $\overline L^{\roman{a}}_{\{\cdot,\cdot\}^{\roman{red}}}$
and hence induces
an
$(A^{\roman{red}},
\overline L^{\roman{a}}_{\{\cdot,\cdot\}^{\roman{red}}})$-module
structure
$$
\overline L^{\roman{a}}_{\{\cdot,\cdot\}^{\roman{red}}}
\longrightarrow
\roman{End}_{\Bobb R}(M^{\roman{red}})
\tag2.8
$$
on
$M^{\roman{red}}$
which is, indeed,
that of an algebraic prequantum module
for $(A^{\roman{red}},\{\cdot,\cdot\}^{\roman{red}})$.
Since the formal differentials which are zero at every point of 
$N^{\roman{red}}$
act trivially on
$M^{\roman{red}}$,
(2.8) factors through a geometric
prequantum module structure
$$
\chi^{\roman{red}} \colon
\overline L_{\{\cdot,\cdot\}^{\roman{red}}}
\longrightarrow
\roman{End}_{\Bobb R}(M^{\roman{red}}).
\tag2.9
$$
\smallskip
For each stratum $N^{\roman{red}}_{(H)}$, restriction of
the reduced V-line bundle $\zeta^{\roman{red}}$
to that stratum
yields a 
V-line bundle  $\zeta_{(H)}$;
since the stratum arises by ordinary symplectic reduction
\cite\sjamlerm,
the space of smooth sections
$M_{(H)} = \Gamma^{\infty}\left(\zeta_{(H)}\right)$ 
of $\zeta_{(H)}$
inherits a
prequantum module structure
$$
\chi_{(H)}
\colon 
\overline L_{\{\cdot,\cdot\}_{(H)}}
\longrightarrow
\roman{End}_{\Bobb R}(M_{(H)})
$$
for
$(A_{(H)},\{\cdot,\cdot\}_{(H)})$.
By naturality, 
the restriction mapping is compatible
with these structures, that is,
the restriction mappings and 
the adjoints 
$\chi^{\sharp}_{\roman{red}}$
and
$\chi^{\sharp}_{(H)}$
of the corresponding prequantum module structures
make 
the requisite diagram
$$
\CD
\phantom{\,}\overline L_{\{\cdot,\cdot\}^{\roman{red}}}
\otimes M^{\roman{red}}\phantom{\{,\}}
@>{\chi^{\sharp}_{\roman{red}}}>>
\phantom{\roman{red}}M^{\roman{red}}
\\
@VVV
@VVV
\\
\phantom{\,\,}\overline L_{\{\cdot,\cdot\}_{(H)}} \otimes M_{(H)}\phantom{\{,\}}
@>{\chi^{\sharp}_{(H)}}>>
\phantom{(H)}M_{(H)} 
\endCD
$$
commutative.
In other words,
the data determine
a stratified prequantum module structure
$(\chi^{\roman{red}}, \{\chi_{(H)}\})$
for the stratified symplectic space
$(N^{\roman{red}},C^{\infty}(N^{\roman{red}}),\{\cdot,\cdot\}^{\roman{red}})$.
This completes the proof of Theorem 2.1. \qed
\enddemo

\medskip\noindent {\bf 3. Quantization}\smallskip\noindent
Let $(N,C^{\infty}(N),\{\cdot,\cdot\})$
be an arbitrary stratified symplectic space, and
let 
\linebreak
$P
\subseteq
\Omega^1(N,\Bobb C)_{\{\cdot,\cdot\}}$
be a stratified complex polarization
(see Section 2 of \cite\kaehler).
We do not recall the general definition
of a stratified complex polarization but,
for intelligibility, we will briefly explain the special case
of {\it complex analytic stratified K\"ahler polarization\/}.
The results of the present paper
will involve only this special case;
we will build up the general theory in terms of
general stratified complex polarizations, though.
Thus suppose that 
the stratified space $N$, endowed with the smooth structure
$C^{\infty}(N)$, is 
{\it fine\/}, that is,
for an arbitrary locally finite
open covering $\Cal U$ of $N$, there is a partition of unity
subordinate to $\Cal U$, and that $N$
is endowed with a complex analytic structure
such that the following hold:
(i) each stratum is complex analytic, that is,
the stratification of $N$
(as a stratified symplectic space)
is a refinement of the complex analytic stratification;
(ii) each (germ of) holomorphic function
belongs to the smooth structure
$C^{\infty}(N,\Bobb C)$;
and
(iii)
for every pair of (germs of) holomorphic functions $f$ and $h$,
the Poisson bracket $\{f,h\}$ is zero.
Then the 
$C^{\infty}(N,\Bobb C)$-submodule
$P$ of $\Omega^1(N,\Bobb C)_{\{\cdot,\cdot\}}$
generated by differentials of the kind
$udf$ where $f$ is holomorphic and $u$ a bump function
in $C^{\infty}(N,\Bobb C)$
(the notion of bump function being interpreted
in terms of the {\it fineness\/} of the structure)
is a complex analytic stratified K\"ahler polarization
for $N$.
See Theorem 2.5 in \cite\kaehler\ for details.
\smallskip
We shall say that a function
$f \in C^{\infty}(N)$
is {\it compatible with the polarization\/} $P$
or, equivalently,
{\it quantizable in the polarization\/} $P$
provided, for every $\alpha \in P$,
$[df,\alpha] \in P$.
When $N$ is a smooth symplectic manifold
and $\{\cdot,\cdot\}$
its ordinary symplectic Poisson structure,
this notion of compatibility 
with the polarization boils down to the usual notion
of classical observable which is
compatible with 
or, equivalently, quantizable in,
the given complex polarization.
See also p. 217 of \cite\souriau.
\smallskip
The extension
of $(\Bobb C,C^{\infty}(N,\Bobb C))$-algebras
determined by the complexified Poisson 2-form 
$\pi \otimes \Bobb C$
arises from  the corresponding extension (1.5.1) by complexification
and has the form
$$
0
@>>>
C^{\infty}(N,\Bobb C)
@>>>
\overline L^{\Bobb C}_{\{\cdot,\cdot\}}
@>>>
\Omega^1(N,\Bobb C)_{\{\cdot,\cdot\}}
@>>>
0.
\tag3.1
$$
Recall that, cf. (1.2.2),
as 
a $C^{\infty}(N,\Bobb C)$-module,
$\overline L^{\Bobb C}_{\{\cdot,\cdot\}}$
is the direct sum of
$C^{\infty}(N,\Bobb C)$ and
$\Omega^1(N,\Bobb C)_{\{\cdot,\cdot\}}$
whence there is
an obvious section $\kappa$ for (3.1).
Now a function
$f \in C^{\infty}(N)$
is  compatible with the polarization $P$
if and only if, for every $\alpha \in P$,
$$
[(f,df),(0,\alpha)] \in \kappa (P).
$$
More formally, we may proceed as follows:
Let $\Omega^1(N,\Bobb C)_P \subseteq \Omega^1(N,\Bobb C)_{\{\cdot,\cdot\}}$
be the  
Lie subalgebra
which consists of all
$\alpha \in \Omega^1(N,\Bobb C)_{\{\cdot,\cdot\}}$
such that, for every $\beta \in P$,
$[\alpha,\beta] \in \Omega^1(N,\Bobb C)_P$;
in other words,
$\Omega^1(N,\Bobb C)_P$ is the {\it normalizer\/}
(in the sense of Lie algebras)
of $P$ in $\Omega^1(N,\Bobb C)_{\{\cdot,\cdot\}}$.
For intelligibility,
we note that, 
under the adjoint $\pi^{\sharp}$ 
from 
$\Omega^1(N,\Bobb C)$ to
the complexified vector fields $\roman{Vect}(N,\Bobb C)$,
the Lie algebra $\Omega^1(N,\Bobb C)_P$ passes to the 
Lie subalgebra
$\roman{Vect}(N,\Bobb C)_P$
of 
$\roman{Vect}(N,\Bobb C)$
which consists 
of complexified vector fields on $N$
that are compatible with the polarization $P$;
here
$\roman{Vect}(N,\Bobb C)$ refers to
the $(\Bobb C,C^{\infty}(N,\Bobb C))$-Lie algebra of
complexified vector fields on $N$.
\smallskip
Let 
$C^P(N,\Bobb C) \subseteq C^{\infty}(N,\Bobb C)$
be the subalgebra of $P$-invariant elements;
in general, non-trivial functions in
$C^P(N,\Bobb C)$ will exist at most locally,
and we should talk about the {\it sheaf of germs of\/}
$P$-invariant functions.
When $N$ is a complex analytic stratified K\"ahler space,
with stratified K\"ahler polarization $P$,
$C^P(N,\Bobb C)$ amounts to the sheaf of germs of
holomorphic functions.
A special case thereof
is that of a smooth K\"ahler manifold where $P$ arises from the ordinary
holomorphic polarization.
In the general case,
the Lie algebras
$\overline L^{\Bobb C}_P$
and
$\Omega^1(N,\Bobb C)_P$ inherit 
$(\Bobb C,C^P(N,\Bobb C))$-Lie algebra structures in an obvious fashion,
the extension (3.1) restricts to an extension
$$
0
@>>>
C^{\infty}(N,\Bobb C)
@>>>
\overline L^{\Bobb C}_P
@>>>
\Omega^1(N,\Bobb C)_P
@>>>
0
$$
of
$(\Bobb C,C^{\infty}(N,\Bobb C))$-algebras,
and
$\overline L^{\Bobb C}_P$
may be viewed as a
sub $(\Bobb C,C^P(N,\Bobb C))$-Lie algebra
of
$\overline L^{\Bobb C}_{\{\cdot,\cdot\}}$.
For later reference, we spell out the following, whose
proof is straightforward and left to the reader.

\proclaim{Proposition 3.2}
A function
$f \in C^{\infty}(N)$
is quantizable in the polarization $P$
if and only 
$(f,df) \in 
\overline L^{\Bobb C}_{\{\cdot,\cdot\}}$
lies in
$\overline L^{\Bobb C}_P$, that is to say,
if and only if
$f$ lies in the pre-image of
$\overline L^{\Bobb C}_P$
under the canonical map $\iota$ from
$C^{\infty}(N)$ to
$\overline L^{\Bobb C}_{\{\cdot,\cdot\}}$. \qed
\endproclaim

Let $(M,\chi)$ be a stratified prequantum module
for
$(C^{\infty}(N),\{\cdot,\cdot\})$.
The composite of the injection
of $P$ into
$\Omega^1(N,\Bobb C)_{\{\cdot,\cdot\}}$
with $\kappa$,
combined with 
the prequantum module structure
$\chi$,
yields a $(C^{\infty}(N,\Bobb C),P)$-module structure 
on 
$M$; 
we denote the subspace
of invariants by $M^P$,
cf. what is said in our paper \cite\souriau.
Likewise,
on each stratum $Y$, 
we have the subspace of invariants $M_Y^{P_Y}$, and the restriction map
induces a linear map
$M^P \to  M_Y^{P_Y}$ of complex vector spaces.
We refer to the system consisting of
$M^P$ and the  restriction maps to the $M_Y^{P_Y}$'s
as the {\it stratified quantum module\/}
determined by the stratified polarization $P$.
Given a stratified quantum module,
the prequantization formula
(1.2.8) induces a representation
of the elements of 
$C^{\infty}(N)$
which are quantizable
in the stratified polarization $P$
by $\Bobb C$-linear operators on
the stratified quantum module,
and this representation
satisfies the conditions (1.2.6) and (1.2.7).
\smallskip
As a side remark
we note that,
more generally, with respect to the Lie-Rinehart algebra
$(C^{\infty}(N,\Bobb C),P)$, the Lie-Rinehart complex
$(\roman{Alt}_{C^{\infty}(N,\Bobb C)}(P,M),d)$
as well as, for every stratum $Y$ of $N$, 
the Lie-Rinehart complexes
$(\roman{Alt}_{C^{\infty}(Y,\Bobb C)}(P_Y,M_Y),d)$,
determine a
system
consisting of the Lie-Rinehart cohomology groups
$\roman H^*(P,M)$ 
together with the  restriction maps to the 
Lie-Rinehart cohomology groups $\roman H^*(P_Y,M_Y)$. 
The system consisting of
$M^P$ and the  restriction maps to the $M_Y^{P_Y}$'s
explained earlier 
(where $Y$ runs through the strata)
boils down to the corresponding zero'th cohomology groups.
The prequantization formula
(1.2.8) now induces a representation
of the elements of 
$C^{\infty}(N)$
which are quantizable
in the stratified polarization $P$
by $\Bobb C$-linear operators on
$\roman H^*(P,M)$ as well as on the 
$\roman H^*(P_Y,M_Y)$'s,
these representations
satisfy the conditions (1.2.6) and (1.2.7) as well,
and the entire system carries the appropriate
costratified structure.
Since we shall not exploit this kind of costratified
structure in the rest of the paper,
we refrain from spelling out details.
An illustration 
for a situation with a single stratum
where ordinary Hodge cohomology groups come into play
will be given shortly.
\smallskip
In particular, suppose that
$(N,C^{\infty}(N),\{\cdot,\cdot\},P)$
is a {\it complex analytic\/} stratified K\"ahler space,
cf. \cite\kaehler, and
let $(M,\chi)$ be a stratified prequantum module
for $(C^{\infty}(N),\{\cdot,\cdot\})$.
We shall refer to
$(M,\chi)$ as a {\it complex analytic\/} stratified prequantum module
provided 
$M$ is the space of sections of a complex V-line bundle
$\zeta$ on $N$ in such a way that
$P$ endows $\zeta$ 
via $\chi$
with a complex analytic structure,
that is to say,
$M^P$ amounts to the sheaf of germs of holomorphic sections of $\zeta$.
\smallskip

Next we show that stratified K\"ahler quantization is 
compatible with passing to the closure of a stratum.
Thus, let $(N,C^{\infty}(N),\{\cdot,\cdot\},P)$
be a stratified K\"ahler space,
let $Y$ be a stratum of $N$, 
and let $(\overline Y,C^{\infty}(\overline Y),\{\cdot,\cdot\}^{\overline Y})$
be the induced stratified symplectic structure on
the closure $\overline Y$ of $Y$,
cf. {\rm (1.4)}.

\proclaim{Proposition 3.3}
The stratified K\"ahler polarization 
$P\subseteq \Omega^1(N,\Bobb C)_{\{\cdot,\cdot\}}$
induces 
a stratified K\"ahler polarization
$P_{\overline Y} \subseteq
\Omega^1(\overline Y,\Bobb C)_{\{\cdot,\cdot\}^{\overline Y}}$
for 
$(\overline Y,C^{\infty}(\overline Y),\{\cdot,\cdot\}^{\overline Y})$.
When
$P$ is complex analytic, so is
$P_{\overline Y}$.
\endproclaim

\demo{Proof}
Let $I$ be the ideal of functions in 
$C^{\infty}(N)$ which vanish on $Y$ (and hence on $\overline Y$).
The canonical projection from $C^{\infty}(N)$ to $C^{\infty}(\overline Y)$
induces an exact sequence
$$
I/I^2 @>>> C^{\infty}(\overline Y) \otimes_{C^{\infty}(N)} \Omega^1(N,\Bobb C)
@>>>
\Omega^1(\overline Y,\Bobb C)
@>>> 0.
$$
The image
$P_{\overline Y}\subseteq \Omega^1(\overline Y,\Bobb C)$ 
of the
induced $C^{\infty}(\overline Y)$-submodule
$C^{\infty}(\overline Y) \otimes _{C^{\infty}(N)} P$
of
$C^{\infty}(\overline Y) \otimes _{C^{\infty}(N)} \Omega^1(N,\Bobb C)$
is a stratified K\"ahler polarization
for 
$(\overline Y,C^{\infty}(\overline Y),\{\cdot,\cdot\}^{\overline Y})$.
When
$P$ is complex analytic, so is
$P_{\overline Y}$. \qed
\enddemo

Let $(M,\chi)$ be a stratified prequantum module for
$(N,C^{\infty}(N),\{\cdot,\cdot\})$,
for any stratum $Y$ of $N$, let
$$
\chi_{\overline Y}\colon
\overline L_{\{\cdot,\cdot\}^{\overline Y}}
@>>>
\roman{End}_{\Bobb R}(M_{\overline Y})
$$
be the induced stratified prequantum module structure
for 
$(\overline Y,C^{\infty}(\overline Y),\{\cdot,\cdot\}^{\overline Y})$,
cf. (1.5.7),
and let $P_{\overline Y} \subseteq
\Omega^1(\overline Y,\Bobb C)_{\{\cdot,\cdot\}^{\overline Y}}$
be the induced stratified K\"ahler polarization.

\proclaim{Theorem 3.4}
For any stratum $Y$,
the morphism 
of stratified prequantum modules
from $(M,\chi)$ to $(M_{\overline Y},\chi_{\overline Y})$
passes to a morphism
of stratified quantum modules
from 
$(M^P,\chi)$ 
to 
$((M_{\overline Y})^{P_{\overline Y}},\chi_{\overline Y})$.
In particular,
for every
pair of strata $Y,Y'$ such that
$Y' \subseteq \overline Y$,
the induced morphism 
of stratified prequantum modules
from $(M_{\overline Y},\chi_{\overline Y})$ to 
$(M_{\overline {Y'}},\chi_{\overline {Y'}})$
passes to a morphism
of stratified quantum modules
from 
$((M_{\overline Y})^{P_{\overline Y}},\chi_{\overline Y})$
to 
$((M_{\overline {Y'}})^{P_{\overline {Y'}}},\chi_{\overline {Y'}})$.
\endproclaim

\demo{Proof} This is straightforward and left to the reader. \qed
\enddemo

We will refer to the costratified complex vector space
which arises from 
the costratified prequantum space
coming from 
a stratified prequantum module on
a stratified symplectic space
(cf. (1.6) above)
by taking invariants, as in Theorem 3.4, 
with reference to a stratified polarization, as
a {\it costratified quantum space\/}.
The linear maps between the constituents
of a costratified quantum space
are not required to be compatible with Hilbert space
structures, though, whatever these structures may be.

\proclaim{Corollary 3.5}
Stratified K\"ahler quantization 
on a (quantizable) complex analytic stratified
K\"ahler space 
$(N,C^{\infty}(N),\{\cdot,\cdot\},P)$
yields a costratified quantum space,
defined on the category
$\Cal C_N$. \qed
\endproclaim

\smallskip
Finally we will show how stratified
quantum modules 
and hence costratified quantum spaces
arise in mathematical nature:
Consider a smooth quantizable (positive) K\"ahler manifold $N$,
viewed  as a stratified symplectic space with a single stratum,
let 
$F \subseteq \roman T^{\Bobb C}N$ be 
its K\"ahler polarization,
and let $P$ be the
(complex analytic) stratified K\"ahler polarization
(in our sense)
arising
as the pre-image 
in 
$\Omega^1(N,\Bobb C)_{\{\cdot,\cdot\}}$
of the space $\Gamma^{\infty} F$ of smooth sections of 
$F$,
with reference to the 
the induced isomorphism
$$
\pi_{\{\cdot,\cdot\}}^{\sharp}\otimes \Bobb C\colon
\Omega^1(N,\Bobb C)_{\{\cdot,\cdot\}} @>>> \roman {Vect}(N,\Bobb C),
$$
cf. (2.1) in \cite\kaehler.
Let $\zeta \colon E \to N$ 
be a prequantum bundle, and
let
$(M,\chi)$ be the corresponding  smooth prequantum module
$(\Gamma^{\infty}(\zeta),\chi_\nabla)$,
$\nabla$ being the corresponding hermitian connection,
cf.  (1.5.2) above;
we view
$(M,\chi)$
as a stratified prequantum module
(with a single stratum).
The line bundle  
$\zeta$ is well known to inherit a
holomorphic structure
whence 
$(M,\chi)$ is a complex analytic stratified prequantum module: 
the quantum module
$M^P$ equals that of
$F$-polarized sections of $\zeta$
in the usual sense
and these are
precisely
the holomorphic sections of $\zeta$.
In this case, the Lie-Rinehart cohomology groups
$\roman H^*(P,M)$ are just the Hodge cohomology groups of $N$ with values
in the holomorphic line bundle $\zeta$.
More generally,
when $N$ is a complex analytic stratified K\"ahler space,
the Lie-Rinehart cohomology 
$\roman H^*(P,M)$ is related with 
the cohomology of 
the sheaf of germs of $\zeta$. 
In the smooth case, the prequantization formula (1.2.8),
applied to $M^P$ and quantizable observables,
then amounts to geometric quantization in the usual sense;
instead of stratified quantum module
we shall then simply say
{\it quantum module\/}.
\smallskip
Let $G$ be a compact Lie group,
and suppose 
that 
{\rm (i)} 
$G^{\Bobb C}$
acts holomorphically on $N$
in such a way that the restriction to
$G$
is  hamiltonian 
and that
{\rm (ii)}
the K\"ahler structure is $G$-invariant;
let
$\mu \colon N \to \fra g^*$
be a corresponding momentum
mapping. 
Suppose that the $G$-action lifts to an action on $\zeta$
preserving the connection.
We have already observed that, for connected $G$,
given the $G$-action on $N$,
the assumption that
$G$ act on $\zeta$ is redundant and
it will suffice
to replace $G$ by an appropriate covering group if necessary.
The prequantum module
$M$ inherits a $G$-action preserving the polarization $P$
and hence the quantum module $M^P$, that is, the space
of holomorphic sections of $\zeta$,
is a complex representation space for $G$.
The quantum module $M^P$
is the corresponding {\it unreduced\/}
quantum state space, 
except that there is no
Hilbert space structure present yet,
and 
{\it reduction after quantization\/},
for the {\it quantum state spaces\/}, amounts to taking
the space
$(M^P)^G$
of $G$-invariant holomorphic sections.
\smallskip
On the other hand, by Proposition 4.2 of \cite\kaehler, 
the (positive) K\"ahler polarization
induces a (positive) 
complex analytic stratified K\"ahler polarization $P^{\roman{red}}$ 
on the reduced space $N^{\roman{red}}$,
with its stratified symplectic Poisson algebra
$(C^{\infty}(N^{\roman{red}}),\{\cdot,\cdot\}^{\roman{red}})$. 
By Theorem 2.1,
the prequantum module $(M,\chi)$ 
passes to a stratified prequantum module
$$
\left(M^{\roman{red}},\chi^{\roman{red}} \colon
\overline L_{\{\cdot,\cdot\}^{\roman{red}}}
\longrightarrow
\roman{End}_{\Bobb R}(M^{\roman{red}})\right)
$$
for the stratified symplectic space
$(N^{\roman{red}},C^{\infty}(N^{\roman{red}}),\{\cdot,\cdot\}^{\roman{red}})$.
{\it Quantization after reduction\/},
for the {\it quantum state spaces\/}, now amounts to taking
the corresponding
{\it reduced\/} quantum module 
or
{\it reduced quantum state space\/}
$(M^{\roman{red}})^{P^{\roman{red}}}$, that is, the space
of 
$P^{\roman{red}}$-invariants in 
$M^{\roman{red}}= \Gamma^{\infty}(\zeta^{\roman{red}})$.
\smallskip
The projection map
from
$M^G$ to
$M^{\roman{red}}=M^G /(IM)^G$
plainly restricts to a linear map
$$
\rho
\colon
(M^P)^G
@>>>
(M^{\roman{red}})^{P^{\roman{red}}}
$$
of complex vector spaces,
referred to as {\it state space comparison map\/} in the introduction
and, 
as far as
the comparison of 
$G$-invariant unreduced 
and reduced quantum observables is concerned,
the statement that
{\it K\"ahler quantization commutes with reduction\/} amounts to the following.

\proclaim{Theorem 3.6}
Let $N$ be a (positive) K\"ahler manifold with a holomorphic
$G^{\Bobb C}$-action whose restriction to $G$
preserves the K\"ahler structure
and is hamiltonian.
Let $f$ be a  $G$-invariant
smooth 
function on $N$
which is quantizable (i.~e. preserves 
the K\"ahler polarization $P$).
Then its class $[f] \in C^{\infty}(N^{\roman{red}})
(=(C^{\infty}(N))^G /I^G)$
is quantizable  
(i.~e. preserves 
the stratified K\"ahler polarization 
$P^{\roman{red}}$) and,
for every
$h \in (M^P)^G$,
$$
\rho(\widehat f (h))
=
\widehat {[f]} (\rho (h)).
$$
\endproclaim

\demo{Proof}
This is seen by direct comparison
of the requisite formulas
(1.2.8) for the unreduced and reduced cases.
We leave the details to the reader. \qed
\enddemo

We now suppose that $\mu$ is an admissible momentum mapping,
that is, that for every $m \in N$ the path of steepest descent through $m$
is contained in a compact set
\cite\kirwaboo\ (\S 9). For example, a proper momentum mapping is admissible.
Another example of an admissible momentum mapping is 
the standard momentum mapping
of a unitary representation space for a compact Lie group;
see Ex. 2.1 in \cite\sjamatwo.
We recall
from \cite\sjamatwo\ 
that, for admissible $\mu$, 
the reduced
bundle
$\zeta^{\roman{red}}$
inherits a holomorphic structure
in the following fashion:
Let $N^{\roman{ss}}\subseteq N$ be the subspace
of semistable points (the points $m$ of $N$ such that the closure of the
$G^{\Bobb C}$-orbit through $m$ intersects the zero level set
$\mu^{-1}(0)$); the inclusion 
$\mu^{-1}(0) \subseteq N^{\roman{ss}}$
then induces a homeomorphism 
$N^{\roman{red}} \to N^{\roman{ss}}\big/\big/G^{\Bobb C}$
\cite\sjamatwo\ (Theorem 2.3).
Likewise,
when 
$E^{\roman{ss}} \subseteq E$
denotes the pre-image of 
$N^{\roman{ss}}$,
the inclusion
$E|_{\mu^{-1}(0)} \subseteq E^{\roman{ss}}$
induces a homeomorphism 
$E^{\roman{red}} =E|_{\mu^{-1}(0)}\big / G
\to E^{\roman{ss}}\big/\big/G^{\Bobb C}$.
According to \cite\roberone,
the sheaf of germs of $G$-equivariant holomorphic sections
of $\zeta|_{N^{\roman{ss}}}$
is a coherent $\Cal O_{N^{\roman{red}}}$-module.
This entails that the 
germs of $G$-equivariant holomorphic sections
of $\zeta|_{N^{\roman{ss}}}$
endow 
$E^{\roman{red}}$ with a complex analytic structure
in such a way that  
$\zeta^{\roman{red}}\colon E^{\roman{red}} \to N^{\roman{red}}$ 
is complex analytic
\cite\sjamatwo\ (Proposition 2.11).
\smallskip
The statement
\lq\lq {\it K\"ahler quantization commutes with reduction\/}\rq\rq\ 
is then completed by the following
two observations which relate the quantum state spaces.

\proclaim{Proposition 3.7}{\rm [\sjamatwo]}
When $\mu$ is admissible and when $N^{\roman{red}}$ has a top stratum
(i.~e. an open dense stratum),
for example when $\mu$ is proper,
the reduced stratified prequantum module $(M^{\roman{red}},\chi^{\roman{red}})$ 
is complex analytic. 
More precisely:
The inclusion
$\Gamma^{\roman{hol}}(\zeta^{\roman{red}})
\to \Gamma^{\infty}(\zeta^{\roman{red}})=M^{\roman{red}}$
of complex vector spaces
identifies $\Gamma^{\roman{hol}}(\zeta^{\roman{red}})$ with
the space $(M^{\roman{red}})^{P^{\roman{red}}}$ of 
$P^{\roman{red}}$-polarized elements of $M^{\roman{red}}$.
\endproclaim

\demo{Proof}
A $P^{\roman{red}}$-polarized element of
$M^{\roman{red}}= \Gamma^{\infty}(\zeta^{\roman{red}})$
is a continuous section
of
$\zeta^{\roman{red}}$
which is holomorphic on each stratum, 
in particular, on the top stratum.
Since, as a complex analytic space,
the reduced space  $N^{\roman{red}}$
is normal, we conclude that a
$P^{\roman{red}}$-polarized
element of $M^{\roman{red}}$ is indeed a holomorphic section
of $\zeta^{\roman{red}}$. \qed 
\enddemo

\proclaim{Theorem 3.8}{\rm [\sjamatwo]}
Under the circumstances
of Theorem {\rm 3.6},
when the momentum mapping $\mu$
is proper, 
in particular when $N$ is compact,
the map $\rho$ is an isomorphism of complex vector spaces.
\endproclaim

\demo{Proof}
This is an immediate consequence of Theorem 2.15 in \cite\sjamatwo. 
In fact, the map $\rho$ is induced by 
the inclusion
$\Gamma^{\roman{hol}}(\zeta^{\roman{red}})
\to \Gamma^{\infty}(\zeta^{\roman{red}})=M^{\roman{red}}$
and the inclusion $N^{\roman{ss}} \subseteq N$.\qed
\enddemo

A version of 
Theorem 3.8 has been established in (4.15) of \cite\naramtwo;
cf. also \cite\sjamafou\ and the literature there, as well as
\cite\ramadthr\ and \cite\ctelethr\ 
for generalizations to higher dimensional sheaf cohomology
and \cite\meinrtwo\ and the literature there
for a generalization in a different direction.

\smallskip
\noindent
{\smc Remark 3.9.}
The statements of Theorems 3.6 and 3.8 are logically independent;
in particular the statement of Theorem 3.6
makes sense whether or not
$\rho$ is an isomorphism, and its proof does
not rely on $\rho$ being an isomorphism.
This raises the issue whether there may exist
states on the reduced level which do {\it not\/} arise
from (invariant) unreduced states.
In fact, when the state space comparison map is not an isomorphism,
interesting things can happen; see \cite\gotaysix.

\smallskip
\noindent
{\smc Remark 3.10.} We do not make any claim as to whether or not
the state space comparison map is unitary.
In the physical context, this issue is, however, vital.
An extensive discussion of the unitarity issue
for the cotangent bundle case 
or, in physicist's terminology, for the Schr\"odinger representation,
may be found in \cite\gotaythr\ together with a thorough treatment
of various examples 
(for the cotangent bundle case)
where the state space comparison map is unitary as
well as of examples where this is not the case.
In \cite\guistetw\ and the subsequent literature 
referenced in \cite\sjamafou,
the question whether or not the state space comparison
map is unitary is not addressed.
Certain cases where unitarity has been established
(under circumstances essentially different from ours)
may be found in \cite\drivhall\ and \cite\bhallone.
In the example treated in the next section, all Hilbert space structures
will derive from the same Hilbert space and the unitarity issue is thereby
circumvented.

\beginsection  4. Holomorphic nilpotent orbits and singular Fock spaces

For illustration,
we will explore a class of examples
involving holomorphic nilpotent orbits;
these arise from the standard simple  Lie algebras 
of hermitian type \cite\kaehler.
On the unreduced level,
the corresponding quantum phase spaces are variants
of ordinary Fock space; reduction
then carries the underlying classical 
phase space
to the closure
of a holomorphic nilpotent orbit. 
Accordingly, we may view the costratified quantum phase space 
arising from K\"ahler quantization
over
the closure of a holomorphic nilpotent orbit  
(cf. (4.6) below)
as a {\it singular\/} Fock space. 
The representations of compact Lie groups which will show up
below are, of course, entirely classical.
What is {\sl new\/} in our approach is the construction
of {\it representations by K\"ahler quantization
on a K\"ahler space with singularities\/}.
\smallskip\noindent
(4.1) {\sl Fock- and related spaces\/}.
Consider $W=\Bobb C^m$,
with its standard complex and K\"ahler structures,
and let $A = C^{\infty}(W)$.
By means of complex coordinates 
$z_1,\dots,z_m$ for $W$---these 
are linear functions on $W$, i.~e. lie in the complex dual
$W^*$ of $W$---the ordinary smooth symplectic
Poisson structure $\{\cdot,\cdot\}$ on the complexification
$A\otimes \Bobb C = C^{\infty}(W,\Bobb C)$
may be described by the formulas
$$
\{z_j,\overline z_k\} = -2i \delta_{j,k},\quad
1 \leq j,k \leq m.
\tag4.1.1
$$
\smallskip
To explain briefly the K\"ahler quantization on $W$
in our framework,
denote by $\zeta \colon W \times \Bobb C \to W$
the trivial line bundle,
with its standard hermitian structure
$\langle \cdot,\cdot\rangle$ which
assigns the value
$\langle \sigma, \sigma'\rangle 
= \phi \overline \phi'$ 
to
two sections $\sigma = \phi \cdot 1$ and
$\sigma' = \phi' \cdot 1$ of $\zeta$.
Recall that
a 1-form $\vartheta \colon \Omega^1(W) \to A \otimes \Bobb C$
is called a (complex) {\it Poisson potential\/}
(for $\{\cdot,\cdot\}$)
provided
$d \vartheta = \pi_{\{\cdot,\cdot\}}$
in the cochain complex computing Poisson cohomology
\cite{\poiscoho,\,\souriau}.
The 
assignments
$$
\vartheta(dz_j) = \frac 12 z_j,
\quad
\vartheta (d\overline z_j) = \frac 12 \overline z_j,
\quad 1 \leq j \leq m,
\tag4.1.2
$$
yield a
Poisson
potential
$\vartheta$
for the Poisson bracket $\{\cdot,\cdot\}$
on 
$A\otimes \Bobb C$
(p. 220 of \cite\souriau);
taking
$M$ to be the space $\Gamma^{\infty}(\zeta)$
of smooth complex sections
of $\zeta$
so that $M$ is essentially a copy of $C^{\infty}(W,\Bobb C)$,
and letting
$$
\chi(a,\alpha)(h) = 
-i \vartheta(\alpha) h
+
\alpha(h) + i a h,\quad h \in C^{\infty}(W,\Bobb C), a \in A, \alpha \in 
\Omega^1(W)_{\{\cdot,\cdot\}},
\tag4.1.3
$$
we obtain a prequantum module structure
$$
\chi
\colon 
\overline L_{\{\cdot,\cdot\}}
\longrightarrow
\roman{End}_{\Bobb R}(M)
\tag4.1.4
$$
for $(A,\{\cdot,\cdot\})$, cf. (1.5) above.
The prequantum module $(M,\chi)$
is plainly geometric and arises from a symplectic structure,
cf. (1.5.3).
Indeed, in the symplectic language,
the symplectic potential which corresponds to 
$\vartheta$
is the ordinary 1-form
$\frac i 4 \sum(z_j d\overline z_j -\overline z_j dz_j)$
which, when $z_j$ is written as  $q_j+ip_j$ ($1 \leq j \leq m$), has the form
$\frac 12 \sum(q_j dp_j - p_j dq_j)$.
The corresponding operator 
of covariant derivative
amounts to the ordinary
hermitian connection on $\zeta$, that is,
for every
$\alpha \in \overline L_{\{\cdot,\cdot\}}$ 
and every
pair of smooth complex sections $\sigma$ and $\sigma'$ of $\zeta$,
$\langle \chi (\alpha)\sigma,\sigma'\rangle+
\langle \sigma,\chi (\alpha) \sigma'\rangle$ equals
$\alpha \langle \sigma, \sigma'\rangle$,
and
$\chi$ 
is characterized by this property.
In view of Proposition 2.3 of \cite\kaehler,
the holomorphic
polarization $P$ 
(or K\"ahler polarization in our sense, cf. Section 3 above,)
is generated by the holomorphic
differentials $dz_1,\dots,dz_m$ and,
smooth complex sections $\sigma$ being identified with smooth 
complex valued functions $\hpsi$ on $W$, 
the resulting quantum module
$M^P$ consists of
smooth complex functions  $\hpsi$ on $W$ satisfying the requirement
$$
0 = \chi(0,dz_j)(\hpsi) = -i \vartheta(dz_j) \hpsi +  \{z_j,\hpsi\},\quad
1 \leq j \leq m,
$$
that is, 
$$
0 = i \{z_j,\hpsi\}+\vartheta(dz_j) \hpsi = 
2 \frac {\partial \hpsi}{\partial \overline z_j} + \frac 12 z_j \hpsi,
\quad
1 \leq j \leq m.
$$
This implies the standard fact that
$M^P$
consists of functions $\hpsi$ in 
$\bold z$ and $\overline {\bold z}$
which may be written in the form
$$
\hpsi(\bold z,\overline {\bold z}) 
= \phi({\bold z}) \roman e^{-\frac {\bold z \overline {\bold z}} 4},
\tag4.1.5
$$
for an entire holomorphic function $\phi$ in $\bold z = (z_1,\dots,z_m)$,
that is to say,
as a module over the algebra 
of entire holomorphic functions
on $W$,
$M^P$
is free, generated by
the function
$\psi_0$ given by the expression
$$
\psi_0(\bold z,\overline {\bold z}) 
= \roman e^{-\frac {\bold z \overline {\bold z}} 4},
\tag4.1.6
$$
and $\langle \psi_0, \psi_0 \rangle 
= \roman e^{-\frac {\bold z \overline {\bold z}} 2}$.
In the physical interpretation,
$\psi_0$ represents the {\it ground state\/}.
The inner product 
$\psi \cdot\psi'$
of two such elements $\psi$ and $\psi'$ 
of $M^P=\Gamma^{\roman{hol}}(\zeta)$, where 
$
\psi(\bold z,\overline {\bold z}) 
= \phi({\bold z}) \roman e^{-\frac {\bold z \overline {\bold z}} 4}
$ and
$\psi'(\bold z,\overline {\bold z}) 
= \phi'({\bold z}) \roman e^{-\frac {\bold z \overline {\bold z}} 4}$,
is given by the standard formula
$$
\psi \cdot\psi' = \int \phi \overline{\phi'}
\roman e^{-\frac {\bold z \overline {\bold z}} 2} \varepsilon_m,
\quad
\varepsilon_m= \frac {\omega^m}{(2 \pi)^m m!},
\tag4.1.7
$$
where $\omega$ refers to th symplectic form,
and the physical Hilbert space,
the {\it bosonic Fock\/} space $\Cal F$,
is the completion of the 
complex vector space
of square integrable functions
of the kind (4.1.5).
For $k \geq 0$,
we will write $\Cal F_k$ for the subspace of $\Cal F$ which consists of
functions of this kind
having $\phi$ a homogeneous degree $k$ polynomial.
It is well known 
that, on each
$\Cal F_k$ ($k \geq 0$),
the integral (4.1.7) converges.
For our purposes, the integral (4.1.7) will provide the requisite
Hilbert space structures. 

\smallskip\noindent
(4.2) {\sl Unreduced observables\/}.
A classical observable $f$,
that is, a function on $W$,
is directly quantizable
(in the holomorphic polarization $P$) provided
$\{z_k,\{z_j,f\}\}$ vanishes for every $1 \leq j,k\leq m$
(p.~219 of \cite\souriau).
In particular, every
classical observable $f$ 
that is at most linear in the $\overline z_j$'s
is quantizable.
The quantization $f \mapsto \widehat f$ of a quantizable
classical observable $f$ is then given by the formula
(1.2.8) which, rewritten in terms of the Poisson potential
$\vartheta$ given by (4.1.2), amounts to
$$
\widehat f(\hpsi) = -i\{f,\hpsi\} + (f- \vartheta (df))\hpsi,
\tag4.2.1
$$
cf. \cite\souriau\ (2.6.7). 
For example, the energy function
$
f_E(\bold z,\overline {\bold z})
=
\frac 12\bold z \overline {\bold z}
$
is quantizable in this polarization
and satisfies
$ \vartheta(df_E) = f_E$;
thus, in view of (4.1.1), its quantization  is
given by
$$
\widehat {f_E}(\hpsi) = -i\{f_E,\hpsi\} 
= - \frac i2 \sum \left(z_j \{\overline z_j, \hpsi\}
+
\overline z_j \{z_j, \hpsi\}\right)
= \sum \left(z_j \frac{\partial \hpsi}{\partial z_j}
-\overline z_j \frac{\partial \hpsi}{\partial \overline z_j}
\right)
\tag4.2.2
$$
whence,
with the notation
$E= \sum z_j \frac{\partial }{\partial z_j}$
for
the ordinary {\smc Euler\/} operator,
for $\hpsi= \phi \psi_0$ in $M^P$, 
cf. (4.1.5),
we have
$$
\widehat {f_E} (\hpsi)
=
E(\phi) \psi_0.
\tag4.2.3
$$
Thus,
when $\phi$ is a homogeneous  degree $k$ 
polynomial, 
$\widehat {f_E} (\hpsi) = k\hpsi$,
the operator $\widehat {f_E}$ has 
the non-negative integers
as its spectrum
and, for $k \geq 0$,
$\Cal F_k$ is the eigenspace 
associated to $k$.
This is of course known to be physically incorrect,
the requisite additional term
arising from the metaplectic correction.
We will address this issue elsewhere.
\smallskip\noindent
(4.3) {\sl Symmetries\/}.
The momentum mapping
(having the value zero at the origin of $W$)
for the 
action of the maximal compact subgroup
$\roman U(m)$ of $\roman{Sp}(m,\Bobb R)\cong \roman{Sp}(W)$
($\roman U(m)$ being the maximal compact subgroup
which fixes the complex structure of $W$)
is well known to be given by
$$
\mu \colon W @>>> \fra u(m)^*,
\quad
\mu^X(\bold z) = \frac i2 \sum x_{j,k}\overline z_j  z_k,
\quad X = [x_{j,k}] \in \fra u(m).
\tag4.3.1
$$
In particular,
for $X =-i\,\roman{Id} \in \fra u(m)$,
$\mu^X$ equals the energy function $f_E$.
For general $X \in \fra u(m)$,
since
$\mu^X$
involves
the $\overline z_j$'s
only linearly,
the function $f=\mu^X$,
viewed as an infinitesimal symmetry,
is quantizable and satisfies
$ \vartheta(df) = f$.
The formula (1.2.9) 
(which refers to infinitesimal symmetries)
then comes down to
$$
\widetilde f(\hpsi) = \{f,\hpsi\} 
\tag4.3.2
$$
and yields the standard $\fra u(m)$-representation 
on 
$\Bobb C[W]=S_{\Bobb C}[W^*]$ ($ = \Bobb C[z_1,\dots,z_m]$)
(by skew-symmetric operators);
this representation integrates to the standard 
$\roman U(m)$-representation
on $S_{\Bobb C}[W^*]$.
\smallskip
Let $\KK$ be a closed subgroup of 
$\roman U(m)$; thus $H$ is a compact subgroup of
$\roman{Sp}(W) = \roman{Sp}(m,\Bobb R)$, and
restricting the 
$\roman U(m)$-representation
yields a 
representation of $\KK$
on 
each $\Cal F_k$ and hence on
$\Cal F$.
Let $G$ be a subgroup of
$\roman{Sp}(W)$
such that $G$ and $\KK$ constitute 
a {\it real\/} ({\it reductive\/}) {\it dual pair\/}
in $\roman{Sp}(W)$ \cite\howeone, and let $\fra g$ be the Lie algebra of $G$,
realized  as a subalgebra of $\fra{sp}(W)$ in the obvious fashion.
We view 
$W$ as the {\it unreduced phase space\/} 
of a classical  system with symmetries given by the
representation of $\KK$ on $W$.
Elements of the Lie algebra $\fra{sp}(W)$
may then be viewed as classical observables
and the elements of $\fra g$
as classical $\KK$-invariant observables.
\smallskip\noindent
(4.4) {\sl Reduction after quantization\/}.
The elements of 
$\fra k = \roman{Lie}(K)$ ($=\fra g \cap \fra u(m)$)
may be viewed as quantizable classical observables.
Restriction yields representations
of $\KK$ and $\fra k$ on 
$\Cal F$ 
and on each 
$\Cal F_k$ ($k \geq 0$)
in such a way that the $\KK$- and $\fra k$-representations
centralize each other;
here and below the Fock space
$\Cal F$
and its homogeneous components depend on the parameter $s$
but we do not indicate this dependence in notation.
We view the $\KK$-representation as a symmetry
and the $\fra k$-representation as
a quantization of classical observables.
The corresponding formulas for
$\fra \kk$ and $\fra k$
then 
result from (1.2.8) and (1.2.9) and, cf. (4.2.2), 
have the form
$$
\widehat f(\hpsi) = -i\{f,\hpsi\}\ 
\text {for}\ \fra k,
\qquad
\widetilde f(\hpsi) = \{f,\hpsi\}\
 \text {for}\ \fra \kk.
\tag4.4.1
$$
{\it Reduction after quantization\/} then 
amounts to passing to the
$\fra k$-representation
on the space $\Cal F^{\KK}$ of
$\KK$-invariants
given by (4.4.1);
this
$\fra k$-representation
decomposes into
representations on the homogeneous degree $k$ constituents
$\Cal F_k^{\KK}$.
\smallskip\noindent
(4.5) {\sl Quantization after reduction\/}.
Let $\mu_{\KK} \colon W \to \fra {\kk}^*$
denote the 
$\KK$-momentum mapping 
having the value zero at the origin;
this momentum mapping is given by the composite of (4.3.1)
with the projection onto
$\fra {\kk}^*$.
By \cite\kaehler\ (Proposition 4.2),
the reduced space $W^{\roman{red}}$
inherits a normal (complex analytic stratified) K\"ahler structure.
\smallskip
We now recall the
following three basic pairs
where we write 
$K = G \cap \roman U(m)$
and where the symmetric constituent 
$\fra p$ of the Cartan decomposition
$\fra g = \fra k \oplus \fra p$
is spelled out explicitly,
for later reference:
\newline\noindent
(4.5.1) 
$(G,\KK,K) = (\roman{Sp}(\ell,\Bobb R),\roman O(s,\Bobb R),\roman U(\ell)),
\ m = 2s\ell, s 
\leq \ell$,
$\fra p \cong S^2_{\Bobb C}[\Bobb C^{\ell}]$;
\newline\noindent
(4.5.2) 
$(G,\KK,K) = (\roman U(p,q),\roman U(s), \roman U(p) \times \roman U(q)),
\ m = s(p+q), s \leq q \leq p$,
$\fra p \cong \roman M_{q,p}(\Bobb C)$;
\newline\noindent
(4.5.3) $(G,\KK,K) =  (\roman O^*(2n),\roman {Sp}(s),\roman U(n)),
\ m = 4ns, s \leq [\frac n2]$,
$\fra p \cong \Lambda^2_{\Bobb C}[\Bobb C^n] \cong \fra o(n,\Bobb C)$;
here $\roman {Sp}(s) = \roman U(s,\Bobb H)$, the unitary group over the
quaternions $\Bobb H$.
\newline\noindent
See Section 5 of \cite\kaehler\ for details and notation.
The pairs (4.5.1),(4.5.2),(4.5.3) correspond  precisely to, respectively,
(5.2)(i),(iii),(iv) in \cite\howeone.
We note that the positive integer
$\ell$, $q$, $[\frac n2]$,
is the real rank of, respectively,
$\roman{Sp}(\ell,\Bobb R)$,
$\roman U(p,q)$,
$\roman O^*(2n)$,
and allowing the parameter $s$ to exceed the real rank of $G$
will not produce
any new examples below.
In case (4.5.1),
the reduced space may be interpreted 
as the classical phase space
of $\ell$ particles in $\Bobb R^s$
with total angular momentum zero.
\smallskip
The general reductive dual pair $(G,H)$ with
$H$ compact arises from taking products of finitely many copies
of the basic pairs.
To simplify the exposition, we will now assume
that $(G,H)$ is any of the three basic pairs.
By Theorem 5.3 of \cite\kaehler,
when
$\fra g^*$ is identified with $\fra g$ via the half-trace pairing
(for the moment any positive multiple of the Killing form would do),
the $G$-momentum mapping $\mu_G \colon W \to \fra g^*$
induces an embedding of
the $\KK$-reduced space 
$W^{\roman{red}}= \mu_{\KK}^{-1}(0)\big /\KK$
into 
$\fra g^*$, and this embedding yields a
normal K\"ahler space
isomorphism
from
$W^{\roman{red}}$
onto the normal K\"ahler space
$
(\overline {\Cal O_s}, C^{\infty}(\overline {\Cal O_s}), 
\{\cdot,\cdot\},P_s^{\roman {red}})
$
whose underlying space
is the closure 
of the holomorphic nilpotent orbit
$\Cal O_s$ in $\fra g$; 
here $P_s^{\roman {red}}$ denotes the (complex analytic)
stratified K\"ahler polarization on
$\overline {\Cal O_s}$ explained in \cite\kaehler.
\smallskip
The holomorphic nilpotent orbits
$\Cal O_0,\dots, \Cal O_r$ are linearly ordered
in such a way that
$$
\{0\}=\Cal O_0 \subseteq \overline{\Cal O_1} 
\subseteq \ldots \subseteq \overline{\Cal O_r}
\tag4.5.4
$$
\cite\kaehler\ (3.3.10);
here $r$ denotes the real rank of $\fra g$.
The orbit $\Cal O_1$ is the
{\it minimal\/} nilpotent orbit (for $G$) which,
in the literature, plays a major role.
The top orbit $\Cal O_r$ is referred to as the 
{\it principal holomorphic\/} nilpotent orbit in \cite\kaehler.
To explain briefly the stratification and the
complex analytic structures recall that, 
with reference to the Cartan decomposition $\fra g = \fra k \oplus \fra p$, 
the symmetric constituent $\fra p$ inherits a complex structure
such that, after complexification so that
$\fra p_{\Bobb C} = \fra p_+ \oplus \fra p_-$
where 
$\fra p_+$ refers to the holomorphic constituent,
the canonical map from $\fra p$ to
$\fra p_+$ is an isomorphism of complex vector spaces;
the orthogonal projection to $\fra p$, restricted to
$\overline{\Cal O_r}$, is a homeomorphism, and in this way
$\overline{\Cal O_r}$ is 
endowed with an
affine complex structure. 
The $G$-momentum mapping $\mu_G \colon W \to \fra g^*$,
combined with the isomorphism onto $\fra g$
induced by the half-trace pairing
followed by the projection to $\fra p$,
amounts to the (complex) Hilbert map of invariant theory---it is
at this stage that the significance of the choice
of half-trace pairing on $\fra g$ comes to the fore;
in this fashion,
the induced map
$\Bobb C[\fra p] \to \Bobb C[W] = S_{\Bobb C}[W^*]$
yields an isomorphism from
the complex affine coordinate ring 
$\Bobb C[\fra p]$ of $\fra p$ onto the $H$-invariants 
$S_{\Bobb C}[W^*]^H$.
This identifies the complex affine coordinate ring 
$\Bobb C [\overline{\Cal O_r}]$
of
$\overline{\Cal O_r}$
with 
$S_{\Bobb C}[W^*]^H$.
For $1 \leq s <r$,
complex analytically, the strata $\Cal O_s$
are the $K^{\Bobb C}$-orbits of
$\fra p \cong \overline{\Cal O_r}$,
and $\overline{\Cal O_s}$ is the corresponding affine subvariety of
$\overline{\Cal O_r} \cong \fra p$.
\smallskip
The algebra
$C^{\infty}(\overline {\Cal O_s})$
is that of Whitney smooth functions on 
$\overline {\Cal O_s}$,
with reference to the embedding into 
$\fra g^*$,
and is {\it not\/} an algebra of ordinary smooth functions;
the Poisson structure comes from this embedding as well.
By Theorem 2.1,
the reduction procedure described in Section 2,
applied to $\zeta$, yields
a stratified prequantum module
$(M^{\roman{red}},\chi^{\roman{red}})$
for 
$\left(C^{\infty}(\overline {\Cal O_s}),\{\cdot,\cdot\}\right)$;
we write 
this stratified prequantum module as
$(M_s^{\roman{red}},\chi_s^{\roman{red}})$.
In the case at hand,
$\zeta$
descends to an ordinary complex line bundle
$\zeta_s=\zeta_s^{\roman{red}}$ 
on $\overline {\Cal O_s}$
which is, in fact, still trivial,
and
$M_s^{\roman{red}}$ is just a free
$C^{\infty}(\overline {\Cal O_s})$-module
of rank 1.
\smallskip
The function
$\psi_0$ 
(cf. (4.1.6))
is invariant under 
$\roman U(m)$ and hence under
$\KK$ 
and
the quantum module
$(M_s^{\roman {red}})^{P_s^{\roman {red}}}$
is
the free
$\Bobb C[\overline {\Cal O_s}]$-module
generated by the (class of the) function $\psi_0$.
Even though the momentum mapping is not proper,
the statement of Theorem 3.7
is still true, that is,
the Hilbert map of invariant theory
induces an isomorphism from
$\Cal F^H = \Bobb C[W^*]^H\langle \psi_0\rangle$
onto $(M_s^{\roman {red}})^{P_s^{\roman {red}}}
= \Bobb C[\overline {\Cal O_s}]\langle \psi_0\rangle$.
In physics, $\psi_0$ amounts to the reduced
ground state.
Furthermore, 
by construction, the formulas (4.4.1) descend; 
hence, 
for every quantizable
$f \in 
C^{\infty}(\overline {\Cal O_s})$
and every
$\psi \in (M_s^{\roman {red}})^{P_s^{\roman {red}}}$,
(1.2.8) and (1.2.9) amount to
$$
\aligned
\widehat f(\hpsi) &= -i\{f,\hpsi\}^{\roman{red}} 
\ \text{(observables)}, 
\\
\widetilde f(\hpsi) &= \{f,\hpsi\}^{\roman{red}}
\ \text{(infinitesimal symmetries)}. 
\endaligned
\tag4.5.5
$$
{\sl These descriptions of the prequantum- and quantum modules
involve only the reduced data and 
make no reference to the
unreduced data.\/}
\smallskip
The 
group $G$ acts on $\overline {\Cal O_s}$ via the adjoint action,
and the subgroup $K$ of $G$ 
is that of transformations preserving the complex analytic stratified
K\"ahler structure.
The inclusion
$\overline {\Cal O_s} \subseteq \fra g$,
combined with the isomorphism $\fra g \cong \fra g^*$
induced by the half-trace pairing,
is a stratified symplectic space momentum mapping
(cf. Section 4 of \cite\kaehler\ for this notion)
for the $G$-action 
on $\overline {\Cal O_s}$
which, combined with the projection to $\fra k^*$,
provides a
stratified symplectic space momentum mapping
$$
\mu \colon \overline {\Cal O_s} @>>> \fra k^*
\tag4.5.6
$$
for the corresponding $K$-action. 
Up to the identification of
$\fra k$ with its dual
(via the half-trace pairing),
this momentum mapping amounts to the
orthogonal projection from 
$\overline {\Cal O_s} \subseteq \fra g = \fra k \oplus \fra p$
to $\fra k$.
Via this momentum mapping,
the elements of the Lie algebra $\fra k$ 
constitute a  Lie subalgebra lying
in the quantizable elements of
$C^{\infty}(\overline {\Cal O_s})$,
and the 
formula (4.5.5)
yields in particular a
$\fra k$-representation
on
$(M_s^{\roman {red}})^{P_s^{\roman {red}}}
= \Bobb C[\overline {\Cal O_s}]\langle \psi_0\rangle$.
\smallskip
The composite of the momentum mapping
(4.5.6)
with 
the infinitesimal generator 
$-Z \in \fra k$
of the central circle subgroup $S^1$
of $K$, viewed
as a linear form on $\fra k^*$,
provides a stratified symplectic  space momentum mapping 
$$
\mu_{-Z} \colon \overline {\Cal O_s} @>>> \roman {Lie}(S^1)^* \cong \bold R
\tag4.5.7
$$
for the hamiltonian $S^1$-action
on $\overline {\Cal O_s}$
obtained from 
letting 
the  circle group $S^1$
act by the inversion map 
from $S^1$ onto
the central circle subgroup 
of $K$; we note that $Z = 2 z$, 
$z$ being the $H$-element of the Lie algebra $(\fra g,z)$
of hermitian type (cf. Section 3 of \cite\kaehler\ for the notion of
$H$-element).
This momentum mapping is
the {\it reduced energy\/}
function
$[f_E]$
(whence the minus sign) and
reveals
certain peculiar features of the reduced system:
The function $[f_E]$
is an element of
$C^{\infty}(\overline {\Cal O_s})$
but
{\it not\/}
an ordinary smooth function,
not even
for $s=r$;
for $s=r$,
$\overline {\Cal O_r} \cong \fra p$
is a complex affine space
and hence a smooth manifold
but the algebra
of functions $C^{\infty}(\overline {\Cal O_r})$
underlying the stratified symplectic structure
is {\it strictly larger\/} than the algebra of ordinary smooth
functions on $\fra p$.
Moreover,
$[f_E]$
is {\it not\/} homogeneous quadratic,
and the reduced energy operator
$\widehat{[f_E]}$ 
(given by the formula (4.5.5))
has only {\it even\/} (non-negative) eigenvalues;
indeed, we can as well compute 
this operator
from the formula
(4.2.3), noticing that
only {\it even\/} entire holomorphic functions
will come into play.
\smallskip\noindent
(4.6) {\sl The costratified quantum space structure\/}.
Let $1 \leq s \leq r$.
Whenever $s'<s$, 
restriction
yields a morphism 
of stratified quantum modules
from
$ (M_s^{\roman {red}})^{P_s^{\roman {red}}}$ to
$(M_{s'}^{\roman {red}})^{P_{s'}^{\roman {red}}}$;
the complex vector spaces 
$M_s^{\roman {red}}$
and
$M_{s'}^{\roman {red}}$
are just the free
$C^{\infty}(\overline {\Cal O_s},\Bobb C)$-
and
$C^{\infty}(\overline {\Cal O_{s'}},\Bobb C)$-modules,
respectively,
in a single generator  whence,
$\overline {\Cal O_{s'}}$ being complex analytically 
an affine subvariety of
$\overline {\Cal O_{s}}$,
the restriction morphism amounts to the
canonical surjection 
$\Bobb C[\overline {\Cal O_s}] \to
\Bobb C[\overline {\Cal O_{s'}}]$
from the affine complex coordinate ring of
$\overline {\Cal O_s}$
to that of
$\overline {\Cal O_{s'}}$.
Thus the resulting costratified quantum space 
for $\overline {\Cal O_s}$ arises from the system
$$
\Bobb C \langle\psi_0 \rangle
@<<<
\Bobb C[\overline {\Cal O_1}]\langle\psi_0 \rangle
@<<<
\ldots @<<< \Bobb C[\overline {\Cal O_s}]\langle\psi_0 \rangle
$$
by Hilbert space completion
where the notation $\psi_0$ for the basis elements is slightly abused.
Here each arrow is actually a morphism of representations
for the corresponding quantizable observables,
in particular,
a morphism of $\fra k$-representations.
Plainly, this structure integrates to
a costratified
$K$-representation, i.~e.
corresponding system of
$K$-representations.
We view the resulting
costratified quantum phase space 
for $\overline {\Cal O_s}$
as a {\it singular\/} Fock space. 
\smallskip\noindent
(4.7) {\sl Quantization commutes with reduction\/}.
By Theorems 3.5
and 3.7,
it makes no difference
whether
we compute the value of an unreduced quantum observable
in an unreduced quantum state
or the value of the corresponding
reduced quantum observable
in the corresponding reduced quantum state.
In particular, this remark applies to
the elements of $\fra k$;
viewed as quantizable reduced observables
via the momentum mapping (4.5.6),
they lie in the reduced Poisson algebra
$C^{\infty}(W^{\roman{red}})\cong
C^{\infty}(\overline {\Cal O_s})$.
Reduction after quantization 
yields the $\fra k$-representation
on the invariants $\Cal F^H=\Bobb C[W]^H\langle \psi_0\rangle$
given by (4.4.1),
quantization after reduction
yields the $\fra k$-representation
on
$\Bobb C[\overline {\Cal O_s}]\langle \psi_0\rangle$
given by (4.5.5),
and the Hilbert map of invariant theory
provides an isomorphism
of $\fra k$-representations
between the two.
\smallskip\noindent
(4.8) {\sl The classical unreduced constant harmonic
oscillator energy phase space\/}.
Let $k \geq 1$ be a positive integer.
The energy function $f_E$ on $W$ 
(given by
$f_E(\bold z,\overline {\bold z})=\frac 12\bold z \overline {\bold z}$)
is the
momentum mapping for the hamiltonian
$S^1$-action on $W$
obtained from 
letting $\alpha \in S^1$
act on $W$
via multiplication by $\alpha^{-1}$.
With respect
to the energy value $k$,  
the reduced space
is a copy of $\Bobb C\roman P^{m-1}$,
endowed with the symplectic form
$k\omega$ where 
$\omega$ 
denotes the symplectic form which is 
the negative of the imaginary part of
the Fubini-Study metric on
$\Bobb C\roman P^{m-1}$.
The $k$'th power $\Cal O(k)=(\Cal O(1))^{\otimes k}$
of the ordinary hyperplane bundle $\Cal O(1)$,
endowed with its hermitian
connection, is a prequantum bundle
for $(\Bobb C\roman P^{m-1}, k \omega)$,
and K\"ahler quantization
yields the (finite dimensional) quantum phase space
$\roman S^k_{\Bobb C}[W^*]=\Cal F_k$.
Indeed,
the line bundle $\Cal O(k)$ actually arises as a special case
of the reduction procedure
for (ordinary smooth) prequantum bundles
spelled out in Theorem 2.1.
Thus K\"ahler quantization on
$\Bobb C\roman P^{m-1}$
recovers $W=\Bobb C^m$
in the sense that
K\"ahler quantization on
$(\roman P[W],k \omega)$ 
picks out the homogeneous 
degree $k$
component of $S_{\Bobb C}[W^*]$;
cf. p.~96 and p.~190 of \cite\woodhous.
In other words,
the symmetric
algebra
$S_{\Bobb C}[W^*] = \Bobb C[z_1,\dots,z_m]$
(viewed as a graded algebra)
being the homogeneous coordinate ring of
$\Bobb C\roman P^{m-1}$,
this homogeneous coordinate ring
amounts to  $\oplus_{k \geq 0}\Gamma^{\roman{hol}}(\Cal O(k))$. 
Corollary 4.11.2 below will 
spell out a similar relationship in a singular situation.
See also the discussion in \cite\schlione\ for related issues.
\smallskip\noindent
(4.9) {\sl Symmetries\/}.
The induced
$\roman U(m)$-action on
$(\Bobb C\roman P^{m-1}, k \omega)$
($k \geq 1$)
is (well known to be) hamiltonian,
the requisite momentum mapping being induced from (4.3.1),
and K\"ahler quantization 
yields the
homogeneous
degree $k$ constituent
$\Cal F_k =\roman S^k_{\Bobb C}[W^*]$, an irreducible summand
of the 
$\roman U(m)$-representation
on $S_{\Bobb C}[W^*] = \Bobb C[z_1,\dots,z_m]$.
Thus this representation is seen as arising by first reducing
with respect to the 
energy function
and quantizing thereafter and,
as an illustration of the principle that
quantization commutes with reduction,
this representation arises as well as the 
Euler operator
eigenspace 
associated to $k$, cf. (4.2) above, 
the Euler operator being the quantized
harmonic oscillator hamiltonian.
\smallskip\noindent
(4.10) {\sl Reduction after quantization\/}.
Let $k \geq 1$.
Reduction after quantization now amounts to passing
to the resulting $\fra k$-representation on
the $\KK$-invariant subspace $\Cal F^{\KK}_k$.
\smallskip\noindent
(4.11) {\sl Quantization after reduction\/}.
In view of Proposition 4.2 and Theorem 10.1 in \cite\kaehler,
symplectic reduction applied to
$(\Bobb C\roman P^{m-1}, k \omega)$
with reference to the induced $\KK$-action
yields the normal compact complex analytic stratified K\"ahler space
$$
Q_{s,k}=(Q_s,C^{\infty}(Q_s), 
\{\cdot,\cdot\}_k^{\roman{red}}, P_s^{\roman{red}});
$$ 
henceforth, when we wish to indicate all the structure 
including the stratified symplectic Poisson bracket,
we will use the notation $Q_{s,k}$.
In case (4.5.1), 
this space may be interpreted 
as the classical phase space
of $\ell$ harmonic oscillators in $\Bobb R^s$
with total angular momentum zero
and constant energy $k$.
As a complex analytic space, $Q_s$
is a projective variety and,
with an abuse of notation, $P_s^{\roman{red}}$ refers to the resulting
complex analytic stratified K\"ahler polarization
on $Q_s$.
\smallskip
With reference to the Cartan decomposition
$\fra g = \fra k \oplus \fra p$,
cf. (4.5) above,
by Theorem 10.1 in \cite\kaehler, when $r$
refers to the {\it real\/} rank
of $\fra g$,
as a complex analytic space,
$Q_r$ is the complex projective space
$\roman P(\fra p)$
and thus amounts to ordinary complex projective space
$\Bobb C \roman P^d$ where
$d= \frac {\ell(\ell+1)}2 -1$, 
$d= pq -1$,
$d= \frac {n(n-1)}2 -1$,
according to, respectively, the cases
(4.5.1), (4.5.2), (4.5.3).
For $1 \leq s \leq r$,
complex analytically,
$Q_s$ arises from
projectivization of the
closure
$\overline {\Cal O_s}$
of the holomorphic nilpotent orbit
$\Cal O_s$,
the strata of $Q_s$
are the $K^{\Bobb C}$-orbits 
in $Q_s$, and their
closures constitute an ascending chain
$Q_1\subseteq Q_2 \subseteq \dots\subseteq Q_s$
of projective varieties.
However, the stratified symplectic Poisson structure
$(C^{\infty}(Q_r),\{\cdot,\cdot\}_k^{\roman{red}})$ on 
$Q_r$ ($=\Bobb C \roman P^d$)
differs from the standard Poisson structure
coming from (a multiple of) the (Fubini-Study) symplectic 
structure in an essential fashion.
\smallskip
The reduction procedure for
prequantum modules
given in Section 2 above,
applied to the prequantum module $M$
arising from the space
of smooth sections
of $\zeta =\Cal O(k)$,
yields a prequantum module 
$M_{s,k}^{\roman{red}}$
for
$(C^{\infty}(Q_s), \{\cdot,\cdot\}_k^{\roman{red}})$.
As a module over $C^{\infty}(Q_s)$,
$M_{s,k}^{\roman{red}}$
consists of
continuous sections
of the bundle
$\zeta_{s,k}^{\roman{red}}$
on $Q_s$
arising from reduction applied to the $k$'th power
$\Cal O(k)$ of the hyperplane bundle
which are smooth on each stratum,
and
$\zeta_{s,k}^{\roman{red}}$
inherits a holomorphic structure.
Theorem 3.8, or a direct argument, entails that the canonical map
$\rho$ 
(coming into play in Theorems 3.6 and 3.8)
identifies the space of holomorphic sections
of
$\zeta_{s,k}^{\roman{red}}$
with that of $\KK$-invariant
holomorphic sections of
$\Cal O(k)$.
\smallskip
For $1 \leq s < r$,
we write
$\iota_{Q_s} \colon Q_s \to
Q_r =\Bobb C\roman P^d$
for the embedding, cf. Theorem 10.1 in \cite\kaehler.
This embedding
determines a homogeneous coordinate ring
$S[Q_s]$ for $Q_s$.

\proclaim{Proposition 4.11.1}
Suppose that $r \geq 2$ (so that $d \geq 1$),
and let
$1 \leq s \leq r$.
For $k \geq 1$,
reduction carries
the holomorphic line bundle $\Cal O(2k)$ 
over $\Bobb C\roman P^{m-1}$
to the 
holomorphic line 
bundle $\Cal O_{Q_s}(k)= \iota^*_{Q_s}\Cal O_{\Bobb C\roman P^d}(k)$ over 
the reduced space $Q_s$
and $\Cal O(2k-1)$ to 
a sheaf 
(or complex V-line bundle)
having no non-zero holomorphic section.
Thus the space of
polarized elements, that is, that of
$P^{\roman{red}}$-invariant ones,
of $M_{s,2k}^{\roman{red}}$ 
is the space of holomorphic sections 
$\Gamma^{\roman{hol}}(\Cal O_{Q_s}(k))$
of 
$\Cal O_{Q_s}(k)$
and that of
polarized elements of $M_{s,2k-1}^{\roman{red}}$ 
is zero.
\endproclaim

\demo{Proof} In view of the naturality of 
the constructions,
it suffices to 
establish the first statement for the case $s=r$.
By invariant theory,
the space of $\KK$-invariant holomorphic sections
of the line bundle $\Cal O(k)$ over $\Bobb C\roman P^{m-1}$
is zero for $k$ odd
and for $k$ even,
the dimension
of the space of $\KK$-invariant holomorphic sections
of $\Cal O(k)$ 
coincides with the
dimension
of the space of holomorphic sections
of the 
line bundle
$\Cal O_{Q_{r}}(\frac k2)$ 
over
$Q_{r}=\Bobb C\roman P^d$.
This implies the first assertion.
The second one makes explicit the present construction of 
quantum module. \qed
\enddemo

Thus, for $k\geq 1$, the compact
normal K\"ahler space $Q_{s,2k}$
is quantizable,
having 
the space of 
$C^{\infty}(Q_{s,2k})$-sections 
of the
holomorphic line bundle
$\Cal O_{Q_s}(k)$
as its stratified
prequantum module. 
Inspection establishes the following.

\proclaim{Corollary 4.11.2}
For 
$1 \leq s < r$,
the restriction 
homomorphism from
\linebreak
$\Gamma^{\roman{hol}}(Q_{r}, \Cal O_{Q_{r}}(1))$
to $\Gamma^{\roman{hol}}(Q_{s}, \Cal O_{Q_s}(1))$
is an isomorphism
of complex vector spaces.
Consequently 
the canonical map
from
the homogeneous coordinate ring
$S[Q_s]$ of 
$Q_s$
to
$\oplus_{k \geq 0}\Gamma^{\roman{hol}}(\Cal O_{Q_s}(k))$
is an isomorphism. \qed
\endproclaim

The second statement of this Corollary says 
that
K\"ahler quantization on $Q_s$ recovers 
$\overline {\Cal O_s}$.  
It also entails 
the (well known) fact
that $Q_s$ is projectively normal,
cf. Ex. 5.14 on p. 126 of \cite\hartsboo\ 
and the discussion in (10.7) of \cite\kaehler. 
\smallskip
For $k\geq 1$,
the K\"ahler quantization procedure
developed in Section 3,
applied to the 
complex analytic stratified K\"ahler space
$Q_{s,2k}$ ($1 \leq s \leq r$),
yields 
the {\it costratified quantum space\/} 
$$
\Gamma^{\roman{hol}}(\Cal O_{Q_1}(k))
@<<<
\ldots
@<<<
\Gamma^{\roman{hol}}(\Cal O_{Q_s}(k)).
$$
Each vector space
$\Gamma^{\roman{hol}}(\Cal O_{Q_{s'}}(k))$
($1 \leq s' \leq s$)
is a representation space
for the quantizable observables in
$C^{\infty}(Q_{s})$,
in particular,
a $\fra k$-representation,
and each
arrow is a morphism of representations;
these arrows are just restriction maps.
This structure globalizes to
a costratified $K$-representation.

\smallskip\noindent
(4.12) {\sl Quantization commutes with reduction\/}.
For $1 \leq s \leq r$,
the 
vector
space
$\Gamma^{\roman{hol}}(\Cal O_{Q_s}(k))$
($k \geq 1$)
coincides with the  finite dimensional space
$\Cal F_k^{\KK}$ of $\KK$-invariants,
and the $\fra k$-representation on
$\Gamma^{\roman{hol}}(\Cal O_{Q_s}(k))$
coincides with the representation (4.4.2)
of $\fra k$ on 
$\Cal F_k^{\KK}$.
The compactness of
$Q_s$
and hence its singular structure,
as made precise in Theorem 10.1 in \cite\kaehler,
are {\it crucial\/} at this stage.
{\it Had we carried out K\"ahler quantization on the top stratum
of
$Q_s$
only, which is in fact a smooth (non-compact) K\"ahler manifold,
we would have obtained an infinite dimensional
quantum phase space\/}
instead of the finite dimensional vector space
$\Gamma^{\roman{hol}}(\Cal O_{Q_s}(k))$.
Thus forgetting the lower strata amounts to a {\it loss
of information\/} and entails {\it inconsistent results\/}.

\smallskip\noindent
(4.13) {\sl The $K$-symmetries on the closures of holomorphic nilpotent 
orbits and their projectivizations\/}.
Let $(\fra g,z)$ be an arbitrary simple Lie algebra of hermitian type,
let $r$ be its real rank and, as before,
let $\{\Cal O_0, \Cal O_1,\dots,\Cal O_r\}$
the the holomorphic nilpotent orbits, ordered in such a way that
$\{0\}=\Cal O_0 \subseteq \overline{\Cal O_1} 
\subseteq \ldots \subseteq \overline{\Cal O_s}$.
Let $1 \leq s \leq r$ and $k \geq 1$.
The compact normal K\"ahler space 
$Q_{s,k}$
arises as well from
the closure
$\overline {\Cal O_s} \subseteq \fra g$ 
of the holomorphic nilpotent orbit
$\Cal O_s$
of $\fra g$ by
stratified symplectic reduction, 
with reference to the 
momentum mapping (4.5.7) and
energy value $k$.
The $K$-action and
stratified momentum mapping (4.5.6) descend to 
a $K$-action 
on $Q_s$ and stratified momentum mapping
$$
\mu_k \colon Q_{s,k} @>>> \fra k^*.
$$
By construction, 
the image $\mu_k (Q_{s,k})$ lies in the hyperplane
$-Z = k$.
A version of the statement of the {\it Kirillov conjecture\/}, 
but for {\it normal K\"ahler spaces\/}
(rather than smooth K\"ahler manifolds)
now takes the following form:
{\sl Those irreducible $K$-representations
which correspond to the coadjoint orbits in the image 
$\mu_{2k}(Q_{s,2k}) \subseteq \fra k^*$
are precisely the representations which occur in 
the stratified K\"ahler quantization 
$\Gamma^{\roman{hol}}(\Cal O_{Q_s}(k))$
of $Q_{s,2k}$.
Consequently
the irreducible representations
occurring in the projective coordinate ring
$S[Q_s]$ or, equivalently, in the affine coordinate ring
$\Bobb C[\overline {\Cal O_s}]$,
correspond bijectively 
to the integral coadjoint orbits in the image 
$\mu(\overline {\Cal O_s}) \subseteq \fra k^*$
of the momentum mapping\/} (4.5.6).
We will show elsewhere that, in fact, 
under the present circumstances,
the statement of the {\it convexity theorem\/}
obtains, that is,
the intersection $\mu_{2k}(Q_{s,2k}) \cap t_+^*$
with a Weyl chamber $t_+^*$ is a convex polytope
which meets exactly 
the coadjoint orbits 
corresponding to the irreducible $K$-representations
in $\Gamma^{\roman{hol}}(\Cal O_{Q_s}(k))$.
\smallskip
These claims
may be justified
by means of
the following observation which also provides further insight:
With a notation introduced in Section 3 of \cite\kaehler,
the Lie algebra $\fra g$ decomposes as
$\fra g = \fra n^-_r \oplus \fra l_r \oplus \fra n^+_r$ 
where $r$ refers to the real rank of $\fra g$;
cf. \cite\kaehler\ (3.3.4(r)). 
Let $\fra k_0 = \fra k \cap \fra l_r$,
let $L_r \subseteq G$ be the subgroup
having Lie algebra
$\fra l_r$, and let
$K_0 = K \cap L_r$.
The adjoint $K$-action on $\fra g$ 
restricts to an action
of $K_0$ on
$\fra n^-_r$ 
and
$\fra n^+_r$ 
in an obvious fashion.
For example,
when $\fra g = \fra{sp}(\ell,\Bobb R)$,
$r=\ell$,
$\fra l_r \cong \fra {gl}(\ell,\Bobb R)$,
$\fra n^-_r \cong \roman S^2_{\Bobb R}[\Bobb R^\ell] \cong \fra n^+_r$,
the space of real symmetric ($\ell \times \ell$)-matrices,
$G = \roman{Sp}(\ell,\Bobb R)$,
$L_r \cong \roman{GL}(\ell,\Bobb R)$,
$K = \roman U(\ell)$,
$K_0 = \roman O(\ell,\Bobb R)$,
and the action thereof on
$\roman S^2_{\Bobb R}[\Bobb R^\ell]$
is the ordinary one.
In the general case,
the vector bundle
$K\times_{K_0}\fra n^+_r @>>> K\big /K_0$
may be identified with the cotangent bundle 
$\roman T^*(K\big /K_0) \to K\big /K_0$
of the compact homogeneous space $K\big /K_0$, in the following fashion:
The total space 
$\roman T^*(K\big /K_0)$  may be written
as $K\times_{K_0}\fra k_0^{\perp}$
where
$\fra k_0^{\perp} \subseteq \fra k^*$
is the annihilator of $\fra k_0$ in $\fra k^*$.
However, the orthogonal projection from
$\fra g = \fra k \oplus \fra p$ to $\fra k$,
restricted to $\fra n_r^+$,
is an isomorphism from 
$\fra n_r^+$ onto a subspace of $\fra k$
which, under the 
(negative of the)
Killing form,
is the orthogonal complement
of $\fra k_0$.
Now, the assignment 
to $(x,Y)\in K\times \fra n_r^+$
of $\roman {Ad}(x)Y \in \fra g$
induces a
$K$-equivariant map
$\Phi$
from
$\roman T^*(K\big /K_0)$ 
to $\fra g$
whose image is the
union
of all pseudoholomorphic nilpotent orbits in $\fra g$,
and the composite
of $\Phi$ 
with the orthogonal projection from
$\fra g= \fra k \oplus \fra p$ to $\fra k$,
followed by the isomorphism
$\fra k \cong \fra k^*$
induced by 
the half-trace pairing,
is the ordinary $K$-momentum mapping 
for the standard hamiltonian $K$-action on $\roman T^*(K\big /K_0)$.
The restriction 
of $\Phi$ to 
$K\times_{K_0}O_r^+$
where $O_r^+ \subseteq \fra n_r^+$
is a suitable \lq\lq positive definite\rq\rq\ 
part of $\fra n_r^+$
is a diffeomorphism onto the top stratum
$\Cal O_r$ of
$\overline {\Cal O_r}$ and thus exhibits
$\Cal O_r$
as a fiber bundle over
$K\big /K_0$
having
$O_r^+$ 
as its fiber.
Likewise, when
$O_r$ is the space of positive semidefinite elements
of $\fra n_r^+$,
the restriction 
of $\Phi$ to 
$K\times_{K_0}O_r$
is a surjection onto the closure
$\overline {\Cal O_r}$ of
$\Cal O_r$.
For example,
when $\fra g = \fra{sp}(\ell,\Bobb R)$
so that $\fra n_r^+$ is the space of real symmetric
($\ell \times \ell$)-matrices
with the ordinary $\roman O(\ell,\Bobb R)$-action,
$O_r^+$ is that of ordinary positive definite ones.

\smallskip\noindent
{\smc Remark 4.14.} 
As a complex analytic space,
$Q_r$ is the complex projective space
$\roman P[\fra p]$ on the complex vector space $\fra p$ and,
by construction, 
$\fra p$ comes with a $K$-representation.
Similarly as in (4.3) above, 
with $\fra p$ instead of $W$,
a choice of $K$-invariant hermitian form
on $\fra p$
then determines a momentum mapping 
from 
$\fra p$ to $\fra k^*$
and hence,
$\roman P[\fra p]$ being endowed with a multiple of the
Fubini-Study symplectic form,
a momentum mapping
$m \colon \roman P[\fra p] \to \fra k^*$. 
The latter 
is a special case of
certain momentum mappings explored
in  \cite\brionone, \cite\nessone\  and elsewhere. In particular,
the Guillemin-Sternberg convexity
result obtains \cite\guisteth,
and ordinary K\"ahler quantization
on the smooth K\"ahler manifolds
$\roman P[\fra p]$
(when the symplectic structure runs through all
positive multiples of the Fubini-Study symplectic structure)
yields the same
unitary $K$-representations
as our
stratified K\"ahler quantization
on the $Q_{r,2k}$'s ($ k \geq 0)$.
However, our
hamiltonian $K$-space structures
on the $Q_{r,2k}$'s
differ from the ordinary hamiltonian $K$-space structures
on the $\roman P[\fra p]$'s in an essential fashion:
For any $k \geq 1$,
the real structure
$(C^{\infty}(Q_{r,2k}),\{\cdot,\cdot\})$ is a stratified symplectic structure
which 
involves continuous functions which are not necessarily smooth and
has the nice feature that it restricts
to a stratified symplectic structure
on any stratum
$Q_{s,2k}$,
and the momentum mapping $\mu_{r,2k}$
is {\it not\/} the map $m$.
\smallskip\noindent
(4.15) {\sl Explicit descriptions\/}.
For $k \geq 0$,
$\Gamma^{\roman{hol}}(\Cal O_{Q_r}(k)) = S_{\Bobb C}^k [\fra p^*]$
(the space of homogeneous degree $k$ polynomial functions on $\fra p$),
and the decomposition of 
$S_{\Bobb C}^k [\fra p^*]$ into its irreducible
$K$-representations is, of course, well known and classical.
We now reproduce suitable highest weight vectors:
\smallskip\noindent
Case (4.5.1): 
Following the procedure on p. 563 of \cite\howeone\   where 
this is done for $K^{\Bobb C} = \roman{GL}(\ell,\Bobb C)$,
introduce coordinates on $\Bobb C^\ell$. These give rise to coordinates
\linebreak
$\{x_{i,j} = x_{j,i}; \, 1 \leq i,j \leq \ell\}$ on 
$\fra p=S_{\Bobb C}^2 [\Bobb C^{\ell}]$,
and the determinants
$$
\delta_1 = x_{1,1},
\ 
\delta_2 = \left | \matrix x_{1,1}& x_{1,2}\\
                           x_{1,2}& x_{2,2}
                   \endmatrix \right|,
\
\delta_3 = \left | \matrix x_{1,1}& x_{1,2} & x_{1,3}\\
                           x_{1,2}& x_{2,2} & x_{2,3}\\
                           x_{1,3}& x_{2,3} & x_{3,3}\\
                   \endmatrix \right|,
\
\text{etc.}
$$
are highest weight vectors for certain $\roman U(\ell)$-representations.
\smallskip\noindent
Case (4.5.2): 
With the obvious notation $x_{i,j}$ for
coordinates on $\fra p =\roman M_{q,p}(\Bobb C)$,
the determinants
$$
\delta_1 = x_{1,1},
\ 
\delta_2 = \left | \matrix x_{1,1}& x_{1,2}\\
                           x_{2,1}& x_{2,2}
                   \endmatrix \right|,
\
\delta_3 = \left | \matrix x_{1,1}& x_{1,2} & x_{1,3}\\
                           x_{2,1}& x_{2,2} & x_{2,3}\\
                           x_{3,1}& x_{3,2} & x_{3,3}\\
                   \endmatrix \right|,
\
\text{etc.}
$$
are highest weight vectors for certain 
$(\roman U(p)\times \roman U(q))$-representations,
cf. e.~g. \cite\zheloboo\ and
\cite\howeone\ p. 567 (where 
this is explained for the complexified group
$\roman{GL}(p,\Bobb C)\times \roman {GL}(q,\Bobb C)$).
\smallskip\noindent
Case (4.5.3): 
Introduce coordinates $x_1,\dots, x_n$ on $\Bobb C^n$ and 
let
$$
\delta_1 = x_1 \wedge x_2,
\ 
\delta_2 = x_1 \wedge x_2 \wedge x_3 \wedge x_4,
\ 
\delta_3 = x_1 \wedge x_2 \ldots \wedge x_5 \wedge x_6,
\
\text{etc.}
$$
These are highest weight vectors for certain $\roman {U}(n)$-representations.
\smallskip
For $1 \leq s \leq r$ and $k \geq 1$,
the 
$K$-representation
$\Gamma^{\roman{hol}}(\Cal O_{Q_s}(k))$
is the sum of 
the irreducible representations having as
highest weight vectors 
the monomials
$$
\delta_1^{\alpha} \delta_2^{\beta} \ldots \delta_s^{\gamma},
\quad 
\alpha +2 \beta + \dots + s\gamma = k,
$$
and the morphism from $\Gamma^{\roman{hol}}(\Cal O_{Q_s}(k))$
to
$\Gamma^{\roman{hol}}(\Cal O_{Q_{s-1}}(k))$
is an isomorphism on 
the span of
those irreducible representations which do not  involve
$\delta_s$
and has the
span of
the remaining ones as its kernel.
The statement referred to above as a version of the Kirillov conjecture
can now be made more explicit
in the following fashion:
{\sl Those irreducible $K$-representations
which correspond to the coadjoint orbits in the image 
$\mu_{2k}(O_{s'}\setminus O_{s'-1}) \subseteq \fra k^*$
of the stratum
$O_{s'}\setminus O_{s'-1}$
($1 \leq s' \leq s$)
are precisely the
irreducible representations having as
highest weight vectors 
the monomials
$\delta_1^{\alpha} \delta_2^{\beta} \ldots \delta_{s'}^{\gamma}$
($\alpha +2 \beta + \dots + s'\gamma = k)$
involving
$\delta_{s'}$ explicitly, i.~e.
with\/} $\gamma \geq 1$.
\smallskip\noindent
{\smc Remark 4.16.} These observations may be interpreted 
by saying that, 
for the compact hamiltonian $K$-spaces
$(Q_{s,2k},\mu_{2k})$,
K\"ahler quantization commutes with reduction
{\sl even though the underlying spaces have singularities\/}.

\smallskip\noindent
{\smc Remark 4.17.} 
Translating back this information to
the spaces $\overline {\Cal O_s}$ ($s \leq r$),
we conclude:
With reference to the stratified $K$-momentum mapping
(4.5.6), {\sl those irreducible $K$-representations
which correspond to the coadjoint orbits in the image 
$\mu(\Cal O_{s'}) \subseteq \fra k^*$
of the stratum
$\Cal O_{s'}$
($1 \leq s' \leq s$)
of $\overline {\Cal O_s}$
are precisely the
irreducible representations having as
highest weight vectors 
the monomials
$\delta_1^{\alpha} \delta_2^{\beta} \ldots \delta_{s'}^{\gamma}$
involving
$\delta_{s'}$ explicitly, i.~e.
with\/} $\gamma \geq 1$.
Thus for the non-compact hamiltonian $K$-spaces
$\overline {\Cal O_s}$,
K\"ahler quantization commutes with reduction
{\sl even though the underlying spaces have singularities\/}.

\bigskip
\centerline{\smc References}
\medskip
\widestnumber\key{999}

\ref \no  \armcusgo
\by J. M. Arms,  R. Cushman, and M. J. Gotay
\paper  A universal reduction procedure for Hamiltonian group actions
\paperinfo in: The geometry of Hamiltonian systems, T. Ratiu, ed.
\jour MSRI Publ. 
\vol 20
\pages 33--51
\yr 1991
\publ Springer Verlag
\publaddr Berlin $\cdot$ Heidelberg $\cdot$ New York $\cdot$ Tokyo
\endref

\ref \no  \brionone
\by M. Brion
\paper Sur l'image de l'application moment
\paperinfo in: \lq\lq S\'eminaire d'alg\`ebre
Paul Dubreuil et Marie-Paule Malliavin\rq\rq, Paris 1986
\jour Lecture Notes in Math.
\vol 1296
\yr 1987
\pages 177--192
\publ Springer Verlag
\publaddr Berlin $\cdot$ Heidelberg $\cdot$ New York
\endref

\ref \no \drivhall
\by B. Driver and B. Hall
\paper Yang-Mills theory and the Bargmann-Segal transform
\jour Comm. Math. Phys.
\vol 201
\yr 1999
\pages 249--290
\endref

\ref \no \emmroeme
\by C. Emmrich and H. Roemer
\paper Orbifolds as configuration spaces of systems with
gauge symmetries
\jour Comm. Math. Phys.
\vol 129
\yr 1990
\pages 69--94
\endref

\ref \no \gotayone
\by M. J. Gotay
\paper Poisson reduction and quantization for the $n+1$ photon
\jour J. of Math. Phys.
\vol 25
\yr 1984
\pages 2154--2159
\endref

\ref \no \gotaythr
\by M. J. Gotay
\paper Constraints, reduction, and quantization
\jour J. of Math. Phys.
\vol 27
\yr 1986
\pages 2051--2066
\endref

\ref \no \gotaysix
\by M. J. Gotay
\paper Negative energy states in quantum gravity
\jour Class. Quant. Grav.
\vol 3
\yr 1986
\pages 487--491
\endref

\ref \no  \guistetw
\by V. W. Guillemin and S. Sternberg
\paper Geometric quantization and multiplicities of group representations
\jour Invent. Math.
\vol 67
\yr 1982
\pages 515--538
\endref

\ref \no \guisteth
\by V. W. Guillemin and S. Sternberg
\paper Convexity properties of the moment mapping
\jour Invent. Math.
\vol 67
\yr 1982
\pages 491--513
\endref

\ref \no \bhallone
\by  B. C. Hall
\paper Geometric quantization and the generalized Segal-Bargmann transform 
for Lie groups of compact type 
\jour Comm. in Math. Phys.
\vol 226
\yr 2002
\pages 233--268
\finalinfo {\tt quant.ph/0012015}
\endref

\ref \no \hartsboo
\by  R. Hartshorne
\book Algebraic Geometry
\bookinfo Graduate texts in Mathematics
 No. 52
\publ Springer Verlag
\publaddr Berlin $\cdot$ G\"ottingen $\cdot$ Heidelberg
\yr 1977
\endref

\ref \no  \howeone
\by R. Howe
\paper Remarks on classical invariant theory
\jour  Trans. Amer. Math. Soc.
\vol 313
\yr 1989
\pages  539--570
\endref

\ref \no \poiscoho
\by J. Huebschmann
\paper Poisson cohomology and quantization
\jour 
J. f\"ur die reine und angewandte Mathematik
\vol  408 
\yr 1990
\pages 57--113
\endref

\ref \no  \souriau
\by J. Huebschmann
\paper On the quantization of Poisson algebras
\paperinfo Symplectic Geometry and Mathematical Physics,
Actes du colloque en l'honneur de Jean-Marie Souriau,
P. Donato, C. Duval, J. Elhadad, G.M. Tuynman, eds.;
Progress in Mathematics, Vol. 99
\publ Birkh\"auser Verlag
\publaddr Boston $\cdot$ Basel $\cdot$ Berlin
\yr 1991
\pages 204--233
\endref

\ref \no \srni
\by J. Huebschmann
\paper
Poisson geometry of certain
moduli spaces
\paperinfo
Lectures delivered at the \lq\lq 14th Winter School\rq\rq, Srni,
Czeque Republic,
January 1994
\jour Rendiconti del Circolo Matematico di Palermo, Serie II
\vol 39
\yr 1996
\pages 15--35
\endref

\ref \no \claustha
\by J. Huebschmann
\paper On the Poisson geometry of certain moduli spaces
\paperinfo in: Proceedings of an international workshop on
\lq\lq Lie theory and its applications in physics\rq\rq,
Clausthal, 1995,
H. D. Doebner, V. K. Dobrev, J. Hilgert, eds.
\publ World Scientific
\publaddr Singapore $\cdot$
New Jersey $\cdot$
London $\cdot$
Hong Kong 
\pages 89--101
\yr 1996
\endref

\ref \no  \oberwork
\by J. Huebschmann
\paper Singularities and Poisson geometry of certain representation spaces
\paperinfo in: Quantization of Singular Symplectic Quotients,
N. P. Landsman, M. Pflaum, M. Schlichenmaier, eds.,
Workshop, Oberwolfach,
August 1999,
Progress in Mathematics, Vol. 198
\publ Birkh\"auser Verlag
\publaddr Boston $\cdot$ Basel $\cdot$ Berlin
\yr 2001
\pages 119--135
\finalinfo{\tt math.DG/0012184}
\endref

\ref \no \kaehler
\by J. Huebschmann
\paper K\"ahler spaces, nilpotent orbits, and singular reduction
\linebreak
\paperinfo {\tt math.DG/0104213}
\jour Memoirs of the AMS (to appear)
\vol \yr \pages 
\endref

\ref \no \severi
\by J. Huebschmann
\paper Singular Poisson-K\"ahler geometry of Severi varieties and 
their ambient spaces
\paperinfo {\tt math.DG/0206143}
\endref

\ref \no \descent
\by J. Huebschmann
\paper Lie-Rinehart algebras, descent, and quantization
\jour Fields Institute Communications (to appear)
\finalinfo{\tt math.SG/0303016}
\endref

\ref \no \kirwaboo
\by F. Kirwan
\book Cohomology of quotients in symplectic and algebraic geometry
\publ Princeton University Press
\publaddr Princeton, New Jersey
\yr 1984
\endref

\ref \no \kostaone
\by B. Kostant
\paper Quantization and unitary representations
\jour Lecture Notes in Math.
\vol 170
\yr 1970
\pages 87--207
\paperinfo In:
Lectures in Modern Analysis and Applications, III, ed. C. T. Taam
\publ Springer Verlag
\publaddr Berlin $\cdot$ Heidelberg $\cdot$ New York
\endref

\ref \no \kralyvin
\by I. S. Krasil'shchik, V. V. Lychagin, and A. M. Vinogradov
\book Geometry of Jet Spaces and Nonlinear Partial Differential Equations
\bookinfo Advanced Studies in Contemporary Mathematics, vol. 1
\publ Gordon and Breach Science Publishers
\publaddr New York, London, Paris, Montreux, Tokyo
\yr 1986
\endref

\ref \no \meinrtwo
\by E. Meinrenken 
\paper Symplectic surgery and the $\roman {Spin}^{\roman c}$-Dirac operator
\jour Advances in Math.
\vol 134
\yr 1998
\pages 240--277
\endref                      

\ref \no \naramtwo
\by M. S. Narasimhan and T. R. Ramadas
\paper Factorization of generalized theta functions
\jour Inventiones
\vol 114
\yr 1993
\pages 565-623 
\endref

\ref \no  \nessone
\by L. Ness
\paper A stratification of the null cone via the moment map
\jour Amer. J. of Math.
\vol 106
\yr 1984
\pages 1231--1329
\endref

\ref \no \ramadthr
\by T. R. Ramadas
\paper Factorization of generalised theta functions II:
The Verlinde formula
\jour Topology 
\vol 35
\yr 1996
\pages  641--654
\endref 

\ref \no \roberone
\by M. Roberts
\paper A note on coherent $G$-sheaves
\jour  Math. Ann.
\vol 275
\yr 1986
\pages 573--582
\endref

\ref \no \schlione
\by M. Schlichenmaier
\paper Singular projective varieties and quantization
\paperinfo in: Quantization of Singular Symplectic Quotients,
N. P. Landsman, M. Pflaum, M. Schlichenmaier, eds.,
Workshop, Oberwolfach,
August 1999,
Progress in Mathematics, Vol. 198
\publ Birkh\"auser Verlag
\publaddr Boston $\cdot$ Basel $\cdot$ Berlin
\yr 2001
\pages 259--282
\endref

\ref \no \sjamatwo
\by R. Sjamaar
\paper Holomorphic slices, symplectic reduction, and multiplicities of 
representations
\jour Ann. of Math.
\vol 141
\yr 1995
\pages 87--129
\endref

\ref \no \sjamafou
\by R. Sjamaar
\paper Symplectic reduction and Riemann-Roch formulas for multiplicities
\jour Bull. Amer. Math. Soc.
\vol 33
\yr 1996
\pages 327--338
\endref

\ref \no \sjamlerm
\by R. Sjamaar and E. Lerman
\paper Stratified symplectic spaces and reduction
\jour Ann. of Math.
\vol 134
\yr 1991
\pages 375--422
\endref

\ref \no \sniabook
\by J. \'Sniatycki 
\book Geometric quantization and quantum mechanics
\bookinfo Applied Mathematical Sciences
 No.~30
\publ Springer Verlag
\publaddr Berlin $\cdot$ Heidelberg $\cdot$ New York
\yr 1980
\endref

\ref \no \sniatone
\by J. \'Sniatycki 
\paper Constraints and quantization
\paperinfo in: Nonlinear partial differential operators
and quantization procedures,
Clausthal 1981, eds. S.~I.~ Anderson and H.~D.~Doebner
\jour Lecture Notes in Mathematics, No.~1037
\pages 301--334
\publ Springer Verlag
\publaddr Berlin $\cdot$ Heidelberg $\cdot$ New York
\yr 1983
\endref

\ref \no \sniawein
\by J. Sniatycki and A. Weinstein
\paper Reduction and quantization for singular moment mappings
\jour Lett. Math. Phys.
\vol 7
\yr 1983
\pages 155--161
\endref

\ref\no \ctelethr
\by C. Teleman
\paper The quantization conjecture revisited
\jour Ann. of Math.
\vol 152
\yr 2000
\pages 1--43
\endref

\ref \no \woodhous
\by N. M. J. Woodhouse
\book Geometric quantization
\bookinfo Second edition
\publ Clarendon Press
\publaddr Oxford
\yr 1991
\endref

\ref \no \zheloboo
\by D. Zhelobenko
\book Compact Lie groups and their representations
\bookinfo Transl. Math. Mono. No. 40
\publ American Math. Soc.
\publaddr Providence R. I.
\yr 1973
\endref
\enddocument